\documentclass{amsart}
\usepackage{amssymb, graphicx, amsmath, epsfig, changebar}

\theoremstyle{plain}
\newtheorem{theorem}{Theorem}

\newtheorem{lemma}{Lemma}
\newtheorem{proposition}{Proposition}
\newtheorem{remark}{Remark}

\theoremstyle{definition}

\theoremstyle{remark}

\numberwithin{equation}{section}

\begin{document}

\title{Polynomial invariants of links in the projective space}
\author{Maciej Mroczkowski}

\address{Department of Mathematics,
         Uppsala University, 751 06 Uppsala, Sweden}
\email{mroczkow@math.uu.se}  

\begin{abstract}
The Homflypt and Kauffman skein modules of the projective space are computed. Both are free and generated by some infinite set of links. This set may be chosen to be $\{ L_n, n \in \mathbb N \cup \{ 0\}\}$, where $L_n$ is an arbitrary link consisting of $n$ projective lines for $n>0$, and $L_0$ is an affine unknot.
\end{abstract}
\maketitle

\section{introduction}
The celebrated Jones polynomial \cite{J} was generalized shortly after its discovery to the Homfly \cite{H} and Kauffman \cite{K} polynomials, for links in $\mathbb R^3$.
The Jones polynomial was extended from the case of links in $\mathbb R^3$ to $\mathbb RP^3$ by Drobotukhina \cite{D1}.

In this paper, the Homfly and Kauffman polynomials are extended to links in $\mathbb RP^3$. The technique we use is similar to the technique of Lickorish and Millett \cite{LM}. They used heavily the notion of descending diagram for links in $\mathbb R^3$. Here, we use the notion of descending diagram for links in $\mathbb RP^3$, introduced by the author \cite{M}.

I wish to thank Oleg Viro for his help.

\subsection {Background: $\mathbb R^3$ case}
The Homfly and Kauffman polynomials of links in $\mathbb R^3$ are defined by skein relations.

For Homfly polynomial, the skein relations are:

(HI) $x^{-1}\vcenter{\hbox{\epsfig{file=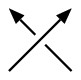}}}-x\vcenter{\hbox{\epsfig{file=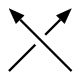}}}=(s-s^{-1})\vcenter{\hbox{\epsfig{file=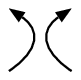}}}$

(HII) $\vcenter{\hbox{\epsfig{file=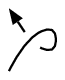}}}=(xv^{-1})\vcenter{\hbox{\epsfig{file=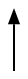}}}$

In these relations, links are presented by their fragments which contain differences from other links under consideration. Moreover, the Homfly polynomial of a link is written simply as the link itself.

The relations (HI) and (HII) together with the assumption that the polynomial is equal to 1 on the unknot $\vcenter{\hbox{\epsfig{file=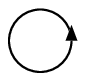}}}$ determine the polynomial for any framed oriented link in $\mathbb R^3$ and this polynomial is invariant under isotopy of such links.

Similarly, for Kauffman polynomial the skein relations are:

(KI) $\vcenter{\hbox{\epsfig{file=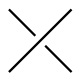}}}+\vcenter{\hbox{\epsfig{file=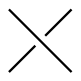}}}=z\left(\vcenter{\hbox{\epsfig{file=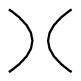}}}+\vcenter{\hbox{\epsfig{file=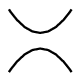}}}\right)$

(KII) $\vcenter{\hbox{\epsfig{file=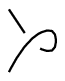}}}=a\vcenter{\hbox{\epsfig{file=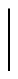}}}$

The relations (KI) and (KII) together with the assumption that the polynomial is equal to $1$ on the unknot $\vcenter{\hbox{\epsfig{file=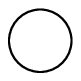}}}$ determine the polynomial for any framed unoriented link in $\mathbb R^3$ and this polynomial is invariant under isotopy of such links.

\subsection{Homflypt and Kauffman skein modules}
Homfly and Kauffman polynomials can be extended to links in any oriented 3-manifold through the notion of skein modules \cite{HP}.

For an oriented 3-manifold $M$, the {\it Homflypt skein module of $M$} is the module over the ring $\mathbb Z[x^{\pm 1},s^{\pm 1},(s-s^{-1})^{-1},v^{\pm 1}]$, generated by isotopy classes of framed oriented links in $M$, with relations (HI) and (HII). The {\it Kauffman skein module of $M$} is the module over $\mathbb Z[a^{\pm 1},z^{\pm 1}]$, generated by isotopy classes of framed unoriented links in $M$, with relations (KI) and (KII).

Note that if $M$ is $\mathbb R^3$, the Homflypt and Kauffman skein modules are free cyclic modules generated by the unknot. The Homfly polynomial (resp. Kauffman polynomial) of a framed oriented link (resp. framed unoriented link) is obtained by expressing the link with the unknot: in the corresponding skein module the link equals to the corresponding polynomial multiplied with the unknot.

\subsection{Main theorems}
\begin{figure}[ht]
\scalebox{1}{\includegraphics{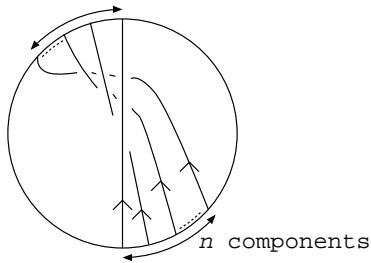}}
\caption{Standard oriented unlink in $\mathbb RP^3$ with $n$ components}
\label{unlink}
\end{figure}

For $n>0$, {\it the standard oriented unlink in $\mathbb RP^3$ with $n$ noncontractible components}, denoted by $L_n$, is the link presented in Figure \ref{unlink} ($\mathbb RP^3$ is represented as a ball $D^3$ with antipodal points of the bounding sphere identified). $L_0$ is the unknot $\vcenter{\hbox{\epsfig{file=circle.eps}}}$. The framing for each $L_n$ (a line framing) is the blackboard one.

\begin{theorem}\label{th_homfly}
The Homflypt skein module of $\mathbb RP^3$ is freely generated by the standard unlinks $L_n$, $n \in \mathbb N\cup \{ 0\}$.
\end{theorem}

This theorem is a consequence of the following:

\begin{theorem}\label{th_homfly_explicit}
To each framed oriented link $L\subset \mathbb RP^3$ a unique element $H(L)\in\mathbb Z[x^{\pm 1},s^{\pm 1},(s-s^{-1})^{-1},v^{\pm 1},z]$ is associated so that $H(L)$ depends only on the isotopy class of $L$, $H(L_n)=z^n$, for $n>0$, $H(L_0)=(v^{-1}-v)/(s-s^{-1})$, and relations (HI) and (HII) hold for $H$.
\end{theorem}

The {\it standard unoriented unlink in $\mathbb RP^3$ with $n$ noncontractible components}, again denoted by $L_n$, is the link presented in Figure \ref{unlink} with the orientations being disregarded. $L_0$ is the unknot $\vcenter{\hbox{\epsfig{file=circle_unoriented.eps}}}$. The framing for each $L_n$ is the blackboard one.

\begin{theorem}\label{th_kauffman}
The Kauffman skein module of $\mathbb RP^3$ is freely generated by the standard unlinks $L_n$, $n \in \mathbb N\cup \{ 0\}$.
\end{theorem}

This theorem is a consequence of the following:

\begin{theorem}\label{th_kauffman_explicit}
To each framed unoriented link $L\subset \mathbb RP^3$ a unique element $K(L)\in\mathbb Z[a^{\pm 1},z^{\pm 1},y]$ is associated so that $K(L)$ depends only on the isotopy class of $L$, $K(L_n)=y^n$, for $n>0$, $K(L_0)=(a+a^{-1}) z^{-1}-1$, and relations (KI) and (KII) hold for $K$.
\end{theorem}

The choices for $H(L_0)$ and $K(L_0)$ are convenient for constructions of $H$ and $K$. These choices are consistent with $H$ and $K$ being equal to 1 for the empty link. Also, with these choices, $H$ and $K$ are multiplicative under disjoint union. Note that in the case of projective links, the disjoint union is well defined on couples consisting of one affine link and one projective link.

\section{Diagrams and basic definitions}
In \cite{M} a notion of descending diagram for links in $\mathbb RP^3$ was introduced.
We use this notion to define Homfly polynomial $H$ (see Theorem \ref{th_homfly_explicit}) for framed oriented links in $\mathbb RP^3$ and Kauffman polynomial $K$ (see Theorem \ref{th_kauffman_explicit}) for framed unoriented links in $\mathbb RP^3$.

\subsection{Diagrams of links in the projective space and nets}
A {\it diagram} of a link in $\mathbb RP^3$ is a disk with a collection of immersed arcs. 
An arc is a compact connected 1-manifold with or without boundary.
The endpoints of arcs with boundary are on the boundary of the disk, divided into pairs of antipodal points and, with this restriction, the arcs are immersed generically. Each double point of the immersions or {\it crossing} of the diagram is endowed with information of over- and undercrossing.

A {\it net} is the projective plane $\mathbb RP^2$ together with a distinguished (projective) line, called {\it the line at infinity}, and a collection of generically immersed circles endowed with information of over- and undercrossing for each double point.
We can map any diagram $D$ of a link to {\it its net}, obtained from $D$ by identifying the antipodal points of the boundary circle of $D$, with the line at infinity being the image of this boundary circle.

If $D$ is a diagram of a link $L$ and $L_b$ is a connected component of $L$, then the projection of $L$ onto $D$ maps $L_b$ onto a collection of arcs. Denote this collection by $b$. We will call $b$ {\it a component} of $D$ (though it may consists of several arcs).
As $H_1 (\mathbb RP^3)=\mathbb Z_2$ there are two types of connected components in a link: 0-homologous and non 0-homologous. 
The corresponding components of a diagram of the link, are said to be {\it 0-homologous} (their images in the net are contractible) and {\it 1-homologous} (their images in the net are not contractible).

If $D$ is a diagram of an oriented link $L$, $D$ gets naturally {\it oriented}: each of its arcs is oriented. If $b$ is a component of $D$ coming from a component $L_b$ of $L$, then the orientation of $L$ gives rise to a cyclic ordering of arcs of $b$ (when one travels on $L_b$ according to the orientation, one meets the arcs of $b$ in this order under the projection of $L$ onto $D$).

\subsection{Arc distance, diagrams descending from $P$ to $Q$}
Let $D$ be a diagram of an oriented link and $b$ a component of $D$. Let $P$ and $Q$ be two points in the interior of some arcs of $b$. Then the {\it arc distance} from $P$ to $Q$ is defined to be the number of times the line at infinity is crossed in the net of $D$, if one travels from the image of $P$ to the image of $Q$ in the net, according to the orientation of the image of $b$ in the net.

Suppose that $X$ is a crossing of $D$ such that at least one of its branches is in $b$. Then the {\it first pass of $X$ from $P$} is, by definition, the branch of $X$ whose image in the net of $D$ is passed first, if one travels from the image of $P$ in the net, according to the orientation of the image of $b$ in the net.

Suppose that $P$ and $Q$ are distinct. $D$ is said to be {\it descending from $P$ to $Q$} if, for every crossing $X$ encountered when traveling from $P$ to $Q$ in the net according to the orientation, the first pass of $X$ from $P$ is an overpass (resp. underpass), if the arc distance from $P$ to this first pass is even (resp. odd). One says that $D$ is {\it descending from $P$ to $P$} if it is descending from $P$ to $Q$, where $Q$ is a point on the same arc as $P$ and such that one can travel on this arc from $Q$ to $P$ according to the orientation without passing any crossing (i.e. $Q$ is just before $P$).

If $D$ is descending from $P$ to $Q$, then, traveling on the net from 
$P$ to $Q$, the encountered arcs are alternatively descending and ascending in the usual meaning for links in $\mathbb R^3$.

\section{Inductive definition of the Homfly polynomial $H$}\label{ind_homfly}
In this section the Homfly polynomial $H$ is constructed on diagrams with a given number of crossings, using the definition of $H$ on diagrams with strictly less crossings.
We will assume that $H$ is constructed on diagrams with less than $n$ crossings and that it has some good properties (see below section \ref{ih}).
The construction of $H$ on diagrams with $n$ crossings requires to endow these diagrams with some extra structure: a basepoint or a couple of basepoints. The construction depends also, in some cases, on an ordering of a set of crossings of a diagram with $n$ crossings.
In the next section we will show that $H$ does not depend on the choice of the extra structure or the choice of ordering, and that it has some good properties on diagrams with $n$ crossings or less.

A {\it standard diagram} of standard oriented unlink $L_n$ is the diagram presented in Figure \ref{unlink}. Let $\mu=(v^{-1}-v)/(s-s^{-1})$.

\subsection{Inductive hypothesis IH$(n-1)$}\label{ih}
There is a function $H$ defined on the set of diagrams with at most $(n-1)$ crossings, taking values in $\mathbb Z[x^{\pm 1},s^{\pm 1},(s-s^{-1})^{-1},v^{\pm 1},z]$ such that:
\begin{enumerate}
\item
$H$ is invariant under those Reidemeister moves that do not increase the number of crossings beyond $n-1$.
\item
$H$ satisfies relations (HI) and (HII).
\item
If $D$ is the standard diagram of standard oriented unlink $L_m$, $m>0$, with at most $(n-1)$ crossings (i.e. $m(m-1)/2\le n-1$), then $H(D)=z^m$.
Also, $H(\vcenter{\hbox{\epsfig{file=circle.eps}}})=\mu$.
\end{enumerate}

\subsection{Diagrams with no crossings}
As the definition of $H$ uses induction on the number of crossings, $H$ is first defined for diagrams with 0 crossings.

Let $D$ be a diagram with 0 crossings. Let $p$ be the number of its 0-homologous components and $m$ the number of its 1-homologous components ($m$ is 0 or 1). Then, by definition:

$(H_1): H(D)=\mu ^p z^m$

For convenience $H$ of the empty link is, by definition, equal to 1 (which agrees with $(H_1)$). Note that $H$ satisfies IH(0). 

\subsection{Diagrams with $n\ge 1$ crossings}
We assume that the inductive hypothesis IH($n-1$) holds true.

The construction of $H$ for diagrams with $n$ crossings is divided into several cases treated in the subsequent subsections.
In each case a diagram $D$ with $n$ crossings is endowed with some extra structure (a basepoint or a couple of basepoints). $D$ together with this structure is denoted by $\dot{D}$. A diagram $\alpha(\dot{D})$ is then defined: it is a diagram obtained from $\dot{D}$ by a series of crossing changes.

The diagram $\alpha(\dot{D})$ has the following property: if $X$ is one of the crossings of $\dot{D}$ that have to be switched to obtain $\alpha(\dot{D})$, and $\dot{D}'$ is the diagram obtained from $\dot{D}$ by switching $X$, then $\alpha(\dot{D})=\alpha(\dot{D}')$.

In the following subsections $H$ is defined on $\alpha(\dot{D})$ for each case (see $(H_3)$ to $(H_6)$).

Suppose that $H$ is already defined on all $\alpha(\dot{D})$. For a based diagram $\dot{D}$, denote by $S(\dot{D})$ the set of crossings of $\dot{D}$ where $\dot{D}$ and $\alpha(\dot{D})$ differ. Let $k$ be the number of elements in $S(\dot{D})$ and $\omega$ a (linear) ordering of $S(\dot{D})$. Denote by $S(\dot{D},\omega)$ the set $S(\dot{D})$ equipped with ordering $\omega$.

We define $H(\dot{D},\omega)$ by induction on $k$. The definition depends on $\omega$.
If $k=0$ then $H$ is already defined.
Otherwise let $\dot{D}'$ be the based diagram obtained from $\dot{D}$ by switching the first crossing in $S(\dot{D},\omega)$ and $D''$ be the diagram obtained by smoothing the same crossing. Let $\omega '$ be an ordering of all crossings of $S(\dot{D'})$ induced by $\omega$. Note that there are $k-1$ elements in $S(\dot{D'})$. $H(\dot{D}',\omega ')$ is defined by induction on $k$ and $H(D'')$ is defined by IH($n-1$).
Now $H(\dot{D},\omega)$ is defined using the relation (HI) on the first crossing in $S(\dot{D},\omega)$ with the help of $H(\dot{D}',\omega ')$ and $H(D'')$. Namely, if $\epsilon$ is the sign of this first crossing, then, by definition:

$(H_2): H(\dot{D},\omega)=x^{2\epsilon} H(\dot{D}',\omega ')+\epsilon x^\epsilon (s-s^{-1}) H(D'')$

\subsection{Simple diagrams}\label{simple_diagrams}
A diagram is said to be {\it simple}, if it has at least one 1-homologous component and any of its crossings involves two different 1-homologous components.
Thus, a simple diagram has no crossings involving 0-homologous components and no self-crossings of 1-homologous components.

A {\it based} simple diagram $\dot{D}$ (the dot indicates that the diagram is based) is a simple diagram $D$ equipped with a couple of {\it basepoints}. The basepoints are antipodal points that lie on the boundary circle of the diagram and that are endpoints of some arc(s) of a 1-homologous component. They are indicated by black dots.

\begin{figure}[ht]
\scalebox{1}{\includegraphics{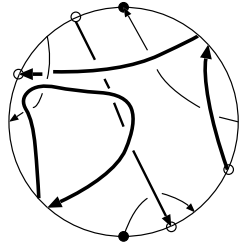}}
\caption{Simple based diagram with 3 components}
\label{simple}
\end{figure}

The antipodal basepoints are also called {\it primary} basepoints. They give rise to an ordering of 1-homologous components and a couple of antipodal basepoints, called {\it secondary}, on each of them (except the component that has the primary basepoints), in the following way: if one travels on the boundary of the disk from the primary basepoints, one encounters, for each 1-homologous component, a couple of endpoints belonging to it. For each such component, the couple encountered for the first time is by definition the couple of secondary basepoints.
It is indicated by a couple of white dots. The 1-homologous components are ordered starting with the component with the primary basepoints, then the component whose secondary basepoints are encountered first and so on, until the component whose secondary basepoints are encountered last.
An example of based simple diagram with a couple of primary basepoints and two couples of secondary basepoints is shown in Figure \ref{simple}.

\begin{figure}[ht]
\scalebox{1}{\includegraphics{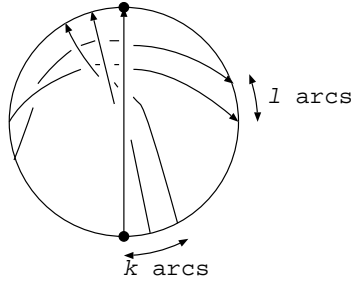}}
\caption{Standard based diagram $\dot{D}_{k,l}$}
\label{unlinkKL}
\end{figure}

Denote by $\dot{D}_{k,l}$ the based diagram shown in Figure \ref{unlinkKL}. It is called {\it standard based diagram}.
A based simple diagram $\dot{D}$ is said to be {\it almost standard}, provided that it can be transformed into some $\dot{D}_{k,l}$ by removing all 0-homologous components and performing a (possibly empty) series of crossing changes.

For any almost standard diagram $\dot{D}$, let $\alpha(\dot{D})$ be the based diagram obtained from $\dot{D}$ by switching all the crossings with sign $-1$. Thus $\alpha(\dot{D})$ has +1 sign at each crossing.

Let $p$ be the number of 0-homologous components of $\dot{D}$ and $m$ the number of its 1-homologous components. Then, by definition:

$(H_3): H(\alpha(\dot{D}))=\mu ^p z^m$

An oriented arc with two endpoints has an {\it initial} endpoint and a {\it final} one. Notice that in an oriented diagram, any couple of antipodal endpoints consists of one initial and one final endpoint.

Let $\dot{D}$ be a simple based diagram.
Consider the endpoint in the couple of primary basepoints of $\dot{D}$ which is initial.
Now travel on the boundary of the disk in the counterclockwise direction from this endpoint and consider each encountered endpoint in the couples of secondary basepoints.
It can be either initial or final. One says that a component is {\it good} if this encountered endpoint is initial and {\it bad} if it is final. The component to which the primary basepoint belongs is always a good component. Note that $\dot{D}_{k,l}$ has $k$ good components followed by $l$ bad ones.

One says that an arc $b$ {\it follows} an arc $a$, if the endpoint of $b$ that is initial, is antipodal to an endpoint of $a$ (which has to be final). An arc is {\it above} (resp. {\it below}) a component if, at each crossing involving the arc and the component, the upper (resp. lower) branch belongs to the arc.

A simple based diagram $\dot{D}$ is said to be {\it descending} if, for any 1-homologous component $b$, the arc of $b$ that contains the basepoint of $b$ which is initial, is above components coming after $b$ according to the order given by the primary basepoints; the arc that follows it is below the same components; and alternating in this way, for all the arcs of $b$.

For any simple based diagram $\dot{D}$, {\it that is not almost standard}, let $\alpha(\dot{D})$ be the based diagram obtained from $\dot{D}$ by crossing changes that make it descending.

Consider the crossings of $\dot{D}_{k,l}$ with sign $-1$. Each such crossing is a crossing between $i$-th good and $j$-th bad component so it may be indexed by $(i,j)$. Let $\omega_{k,l}$ be the lexicographical order of the set of all crossings with sign $-1$, indexed in this way.

Let $p$ be the number of 0-homologous components of $\dot{D}$, $k$ the number of good components and $l$ the number of bad components. Then, by definition:

$(H_4): H(\alpha(\dot{D}))=\mu ^p H(\dot{D}_{k,l},\omega_{k,l})$

Note that in the case of an almost standard diagram $\dot{D}$, $\alpha (\dot{D})$ is not descending unless there are only good components in $\dot{D}$. In what follows, it is more convenient to have $\dot{D}_{k,l}$ descending rather than $\alpha(\dot{D}_{k,l})$.

\subsection{Non simple diagrams}\label{nonsimple_diagrams}
\begin{figure}[ht]
\scalebox{1}{\includegraphics{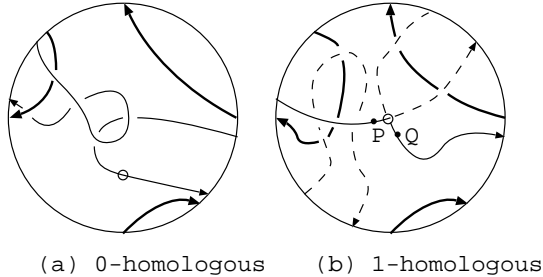}}
\caption{Non simple based diagrams}
\label{nonsimple}
\end{figure}

A notion of based diagram, different from the one used for simple diagrams, is needed in order to construct $H$ for non simple diagrams.
Let $D$ be a non simple diagram. A basepoint can be of two different kinds:
\begin{enumerate}
\item
A point on a 0-homologous component in the interior of some arc and distinct from any crossing {\it if this component is involved in some crossings} (see Figure \ref{nonsimple}(a)).
\item
A self-crossing of a 1-homologous component (see Figure \ref{nonsimple}(b)).
\end{enumerate}

As before, a diagram $D$ equipped with a {\it basepoint} is denoted by $\dot{D}$.

First, we consider $\dot{D}$ for which the basepoint is on a 0-homologous component. $\dot{D}$ is said to be {\it descending}, if it is descending from the basepoint to the basepoint.
If $\dot{D}$ has a basepoint on a 0-homologous component, the diagram $\alpha(\dot{D})$ is the based diagram obtained from $\dot{D}$ by the crossing changes that are necessary to make it descending.

Let $D'$ be the diagram obtained from $\dot{D}$ by removing the component with the basepoint. Let $w$ be the sum of all signs at all self-crossings of the component with the basepoint in $\alpha(\dot{D})$. Then, by definition:

$(H_5): H(\alpha(\dot{D}))=\mu (xv^{-1})^w H(D')$

Finally, we define $H$ in the case when the basepoint of $\dot{D}$ is a self-crossing of a 1-homologous component.

Let $X$ be a self-crossing of a 1-homologous component, say $b$. A smoothing at $X$ according to any orientation of $b$ gives rise to two components: one 0-homologous, the other 1-homologous. The {\it dashed part determined by $X$} is the part of $D$ corresponding to the 0-homologous component.

Let $P$ and $Q$ be two points on $b$ such that, with respect to the orientation, $P$ is just before $X$, $Q$ is just after $X$ and $P$ and $Q$ are not in the dashed part determined by $X$ (see Figure \ref{nonsimple}(b)).
Then $\dot{D}$ is said to be {\it descending}, if it is descending from $P$ to $Q$ (i.e. the dashed part is descending).
If $\dot{D}$ has a basepoint which is a self-crossing of a 1-homologous component, the diagram $\alpha(\dot{D})$ is the based diagram obtained from $\dot{D}$ by the crossing changes which are necessary to make it descending.

Let $D'$ be the diagram obtained from $\dot{D}$ by removing the dashed part determined by the basepoint. Let $w$ be the sum of signs at all self-crossings of the dashed part determined by the basepoint in $\alpha(\dot{D})$, including the basepoint (which is a self-crossing). Then, by definition:

$(H_6): H(\alpha(\dot{D}))=(xv^{-1})^w H(D')$

\section{Independence of $H$ on choices, invariance under Reidemeister moves}
In this section, we show that the definition of $H$ on diagrams with $n$ crossings does not depend on the choices that are involved in it (choice of basepoint and ordering of some crossings). We also prove that the relations (HI) and (HII) are satisfied for diagrams with $n$ crossings and that $H$ is invariant under those Reidemeister moves which do not increase the number of crossings beyond $n$. Thus, assuming that $H$ satisfies the inductive hypothesis IH($n-1$), it is proven that it satisfies IH($n$).

\subsection{Relations (HI) and (HII)}
\begin{lemma}[Independence on ordering]\label{order_independent}
Let $\dot{D}$ be a based diagram with $n$ crossings.
Let $\omega$ and $\omega '$ be two orderings of the set of crossings of $\dot{D}$ that differ between $\dot{D}$ and $\alpha(\dot{D})$.

Then $H(\dot{D},\omega)=H(\dot{D},\omega ')$.
\begin{proof}
By induction on the number of crossing differences between $\dot{D}$ and $\alpha (\dot{D})$ it is sufficient to prove that $H$ does not change if one switches the first two crossings according to $\omega$, say $C_1$ and $C_2$ with respective signs $\epsilon_1$ and $\epsilon_2$.

Denote by $\sigma_1 \dot{D}$ (resp. $\sigma_2 \dot{D}$) the diagram obtained from $\dot{D}$ by switching crossing $C_1$ (resp. $C_2$). Also, denote by $\eta_1 \dot{D}$ (resp. $\eta_2 \dot{D}$) the diagram obtained from $\dot{D}$ by smoothing at crossing $C_1$ (resp. $C_2$). First consider the sequence in which $C_1$ is switched before $C_2$:

\begin{multline*}
H(\dot{D},\omega) = x^{2\epsilon_1}H(\sigma_1 \dot{D})+\epsilon_1 x^{\epsilon_1}(s-s^{-1})H(\eta_1 \dot{D})\\
=x^{2\epsilon_1}x^{2\epsilon_2}H(\sigma_2\sigma_1 \dot{D})+x^{2\epsilon_1}\epsilon_2 x^{\epsilon_2} (s-s^{-1})H(\eta_2\sigma_1 \dot{D})+\epsilon_1 x^{\epsilon_1} (s-s^{-1})H(\eta_1 \dot{D})
\end{multline*}

And switching $C_2$ before $C_1$:

$$H(\dot{D},\omega')=x^{2\epsilon_2}x^{2\epsilon_1}H(\sigma_1\sigma_2 \dot{D})+x^{2\epsilon_2}\epsilon_1 x^{\epsilon_1}(s-s^{-1})H(\eta_1\sigma_2 \dot{D})+\epsilon_2 x^{\epsilon_2}(s-s^{-1})H(\eta_2 \dot{D})$$

The first terms are equal. By inductive hypothesis:

$$H(\eta_1 \dot{D})=x^{2\epsilon_2}H(\sigma_2\eta_1 \dot{D})+\epsilon_2 x^{\epsilon_2}(s-s^{-1})H(\eta_2\eta_1 \dot{D})$$
and
$$H(\eta_2 \dot{D})=x^{2\epsilon_1}H(\sigma_1\eta_2 \dot{D})+\epsilon_1 x^{\epsilon_1}(s-s^{-1})H(\eta_1\eta_2 \dot{D})$$

Substituting these expressions above, one sees that $H(\dot{D},\omega)=H(\dot{D},\omega')$
\end{proof}
\end{lemma}

\begin{lemma}\label{easy_shortening}
Let $D$ be a diagram with at most $n-1$ crossings. Suppose that there is a self-crossing $X$ of a component $a$ of $D$ such that: the arc distance from the upper branch to the lower branch of $X$ is even and $D$ is descending from the upper branch to the lower branch of $X$.

Let $a'$ be the part of $a$ that is covered if one travels in the net from the upper branch to the lower branch of $X$.
Let $D'$ be the diagram obtained from $D$ by erasing $a'$.
Let $w$ be the sum of signs of crossings at which both branches belong to $a'$ (including $X$).

Then $H(D)=(xv^{-1})^w H(D')$.

\begin{proof}
If $a$ is 1-homologous then by IH($n-1$) one can calculate $H(D)$ by putting a basepoint equal to the self-crossing $X$. The conclusion of the lemma follows from $(H_6)$.

Suppose now that $a$ is 0-homologous. The lemma is proven by induction. Let $k$ be the number of crossings in $D$.
For $k=1$ (there has to be at least one crossing in $D$) the lemma is true because of $(H_5)$.
Suppose that the lemma is true for $k<l\le n-1$ and that $D$ has $k=l$ crossings.
The orientation of $D$ induces orientations on the branches of $X$. To compute $H(D)$, put a basepoint on the upper branch of $X$, just before the crossing. $D'$ is obtained from $D$ by removing $a'$, and it inherits the basepoint from $D$.

Now, in computing $H(D)$ and $H(D')$ one can use relation (HI) on crossings that appear both in $D$ and $D'$ since the part $a'$ is already descending with respect to the basepoint in $D$. Let $Y$ be a crossing in $D$ and $D'$ that needs to be changed in order to make the components with basepoints descending. Let $\sigma$ be the operation of switching $Y$, and $\eta$ the operation of smoothing $Y$ with respect to the orientation. For simplicity, assume that the sign of $Y$ is +1. Then, by (HI):

$$x^{-1} H(D)-x H(\sigma D)=(s-s^{-1}) H(\eta D)$$
$$x^{-1} H(D')-x H(\sigma D')=(s-s^{-1}) H(\eta D')$$

By induction, as $\eta D$ has $k-1$ crossings, $H(\eta D)=(xv^{-1})^w H(\eta D')$. Also $H(\sigma D)=(xv^{-1})^w H(\sigma D')$, which is proven easily by induction on the number of crossings that need to be changed in order to make $D$ and $D'$ with basepoints descending.

Thus $H(D)=(xv^{-1})^w H(D')$.
\end{proof}
\end{lemma}

\begin{proposition}[Homfly relations]\label{relation}
The relation (HII) holds for $H$ in the case when the diagram on the left has $n$ crossings.
The relation (HI) holds for $H$ in the case when the two based diagrams on the left have $n$ crossings and these diagrams have the same basepoint(s).
\begin{proof}
For the relation (HII), using the definition of $H$ and IH($n-1$) the statement can be verified easily by considering separately two cases: 
\begin{itemize}
\item
the kink appearing in (HII) is either in a dashed part determined by the basepoint, or in the 0-homologous component to which the basepoint belongs
\item
the kink is not as in the first case
\end{itemize}

The proof for (HI) is more difficult.
Suppose that $\dot{D}$ is a based diagram with $n$ crossings. Let $C$ be a crossing of $\dot{D}$.
Let $\sigma \dot{D}$ be the diagram obtained from $\dot{D}$ by switching $C$, and $\eta \dot{D}$ the diagram obtained from $\dot{D}$ by smoothing $C$. Without lost of generality one may suppose that the sign of $C$ is +1 (otherwise the roles of $\dot{D}$ and $\sigma \dot{D}$ are switched). We want to show that:

$$x^{-1}H(\dot{D})-xH(\sigma \dot{D})=(s-s^{-1})H(\eta \dot{D})$$

First suppose that $\dot{D}$ is simple, or that $C$ is in the dashed part determined by the basepoint, or that $C$ is in the 0-homologous component containing the basepoint.
In these cases, in the definition of $H(\dot{D})$ or $H(\sigma \dot{D})$ the relation (HI) is used at the crossing $C$. By Lemma \ref{order_independent} it can be used at the beginning. Thus (HI) holds for $H$ in these cases.

Suppose now that the crossing $C$ is such that none of its branches is in the 0-homologous component with basepoint of $\dot{D}$ or in the dashed part determined by the basepoint.

If $\dot{D}$ is descending, let $D'$ be the diagram obtained from $\dot{D}$ by erasing the 0-homologous component on which the basepoint lies, or the dashed part determined by the basepoint. Note that the crossing $C$ can be naturally viewed as a crossing in $D'$. As $D'$ has at most $n-1$ crossings, (HI) holds by IH($n-1$) for $H(D')$, $H(\sigma D')$ and $H(\eta D')$.

But $H(\dot{D})$, $H(\sigma \dot{D})$ and $H(\eta \dot{D})$ are expressed respectively with $H(D')$, $H(\sigma D')$ and $H(\eta D')$ in the same way. This is so because of the definition of $H$ except in the case when the basepoint is a self-crossing of a 1-homologous component and $C$ lies on the same component. But in this case, Lemma \ref{easy_shortening} can be used to express $H(\eta \dot{D})$ with $H(\eta D')$.
Thus (HI) holds for $H(\dot{D})$, $H(\sigma \dot{D})$ and $H(\eta \dot{D})$.

Now, if $\dot{D}$ is not descending, (HI) for $H(\dot{D})$, $H(\sigma \dot{D})$ and $H(\eta \dot{D})$ is proved by induction on the number of crossings that have to be changed in order to make it descending. At a crossing that has to be changed, relation (HI) allows to express $H(\dot{D})$, $H(\sigma \dot{D})$ and $H(\eta \dot{D})$, with diagrams with less crossings to be changed, for which (HI) holds by induction, and diagrams with at most $n-1$ crossings for which (HI) holds by IH($n-1$).
\end{proof}
\end{proposition}

The following lemma is a consequence of relation (HI):

\begin{lemma}[basepoints and crossing changes]\label{crossing_independence}
Let $\dot{D}_1$ and $\dot{D}_2$ be two based diagrams with $n$ crossings, that differ only by the position of the basepoint or basepoints.
Let $\dot{D}'_1$ be the diagram obtained from $\dot{D}_1$ by a crossing change and $\dot{D}'_2$ be the diagram obtained from $\dot{D}_2$ by the same crossing change. Suppose that $H(\dot{D_1})=H(\dot{D_2})$.

Then $H(\dot{D}'_1)=H(\dot{D}'_2)$
\begin{proof}
The (HI) relation allows to express $H(\dot{D}_1)$ with $H(\dot{D}'_1)$ and $H$ of a third diagram with $n-1$ crossings. Similarly, it allows to express $H(\dot{D}_2)$ with $H(\dot{D}'_2)$ and $H$ of the same diagram with $n-1$ crossings.
By induction, for the diagram with $n-1$ crossings, $H$ does not depend on the choice of basepoint or basepoints. Thus, if $H(\dot{D}_1)=H(\dot{D}_2)$, then $H(\dot{D}'_1)=H(\dot{D}'_2)$.
\end{proof}
\end{lemma}

\subsection{Basepoints for 0-homologous components}
\begin{lemma}[moving the basepoint]\label{moving_basepoint}
Suppose that $P$ is a basepoint lying on a 0-homologous component of a based diagram $\dot{D}$ with $n$ crossings. Let $c$ be the arc on which $P$ lies.
 
Then $H$ does not change if $P$ is moved on $c$.
\begin{proof}
Denote by $a$ the 0-homologous component with $P$.
It is sufficient to prove that $H$ does not change when the basepoint passes through a crossing as in Figure \ref{h_basepoint}.

\begin{figure}[ht]
\scalebox{1}{\includegraphics{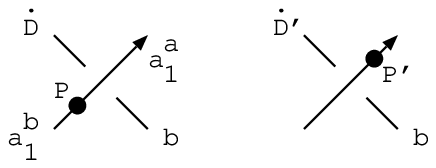}}
\caption{}
\label{h_basepoint}
\end{figure}

Using Lemma \ref{crossing_independence}, one may suppose that the diagram $\dot{D}$ on the left of Figure \ref{h_basepoint} is descending.
Starting from $P$ and traveling on the net of $\dot{D}$ according to the orientation, denote the successive arcs encountered by $a_1, a_2, ..., a_l$. Furthermore denote the part of $a_1$ that comes after $P$ by $a_1^a$ and the remaining part by $a_1^b$. As it was shown in \cite{M} one has the following:

$$a_2\le a_4\le a_6\; ...\le a_1^b\; ...\le a_5\le a_3\le a_1^a$$

where $a_i\le a_j$ means that at each crossing involving $a_i$ and $a_j$ the branch in $a_i$ is under the branch in $a_j$, and $\le$ is transitive.

Now $\dot{D}'$ may be descending or not.
It is not descending if and only if the branch $b$ (see Figure \ref{h_basepoint}) is a part of $a_1^a$, $a_1^b$ or $a_k$ with $k$ odd.
If $\dot{D}'$ is not descending, it becomes descending if one switches the crossing in Figure \ref{h_basepoint}.

Let $\sigma \dot{D}'$ be the based diagram obtained from $\dot{D}'$ by this switching and $\eta \dot{D}'$ the diagram obtained from $\dot{D}'$ by smoothing the same crossing. Note that in $\eta \dot{D}'$, $a$ becomes a link with two 0-homologous components $a^1$ and $a^2$ where $a^1$ contains $P'$ and $a_2$ contains $P$ ($P$ and $P'$ can be naturally viewed in $\eta \dot{D}'$).
Notice that $\eta \dot{D}'$ with basepoint $P'$ is descending.

Denote by $D_a$ the diagram obtained from $\dot{D}'$ by removing $a$. Notice that $D_a$ can be obtained from $\eta \dot{D}'$ by removing $a^1$ and $a^2$. Denote by $D_{a^1}$ the diagram obtained from $\eta \dot{D}'$ by removing $a^1$. Notice that $D_{a^1}$ with basepoint $P$ is descending. Finally, let $\epsilon=\pm 1$ be the sign of the crossing in Figure \ref{h_basepoint}.

For a 0-homologous component $c$ denote by $w(c)$ the sum of signs at all crossings for which both branches are in $c$.

Then:

$H(\dot{D})=\mu (xv^{-1})^{w(a)} H(D_a)$

$H(\eta \dot{D})=\mu (xv^{-1})^{w(a_1)} H(D_{a^1})=\mu^2 (xv^{-1})^{w(a_1)+w(a_2)}H(D_a)$

$H(\dot{D}')=x^{2\epsilon}H(\sigma \dot{D}')+\epsilon x^\epsilon (s-s^{-1})H(\eta \dot{D})=x^{2\epsilon}\mu (xv^{-1})^{w(a)-2\epsilon}H(D_a)+\epsilon x^\epsilon (s-s^{-1})\mu^2 (xv^{-1})^{w(a_1)+w(a_2)}H(D_a)$

As $\mu=(v^{-1}-v)/(s-s^{-1})$:

$H(\dot{D}')=(x^{w(a)}\mu v^{-w(a)+2\epsilon}+\epsilon x^\epsilon \mu(v^{-1}-v)(xv^{-1})^{w(a_1)+w(a_2)})H(D_a)$

As $w(a_1)+w(a_2)=w(a)-\epsilon$:

$H(\dot{D}')=\mu x^{w(a)} v^{-w(a)}(v^{2\epsilon}+\epsilon v^{\epsilon-1}-\epsilon v^{\epsilon+1})H(D_a)=\mu (xv^{-1})^{w(a)}H(D_a)$

Because, whether $\epsilon=1$ or $-1$, one checks easily:

$v^{2\epsilon}+\epsilon v^{\epsilon-1}-\epsilon v^{\epsilon+1}=1$

Thus $H(\dot{D})=H(\dot{D}')$
\end{proof}
\end{lemma}

\subsection{Good Reidemeister moves}
\begin{figure}[ht]
\scalebox{1}{\includegraphics{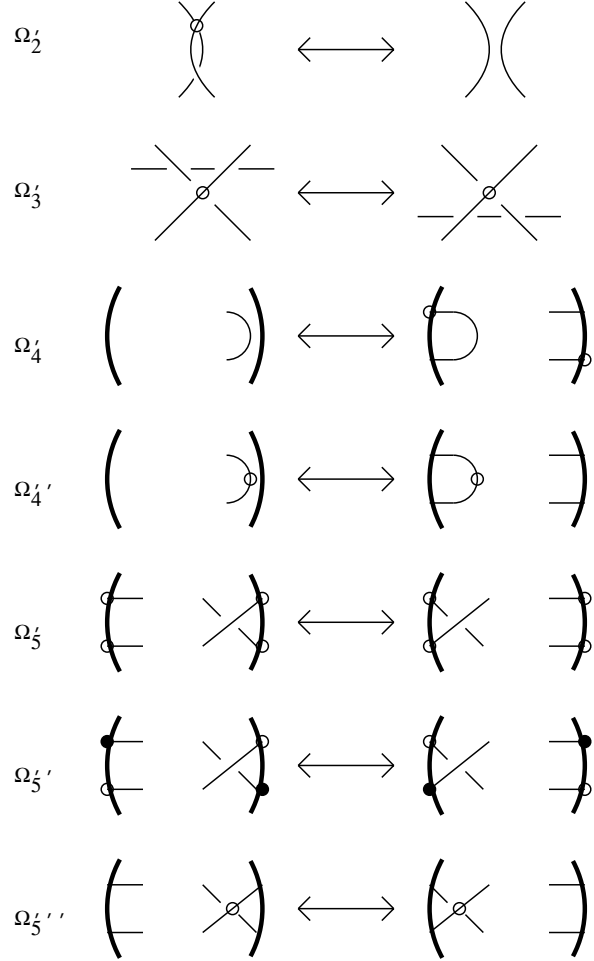}}
\caption{Bad Reidemeister moves}
\label{badmoves}
\end{figure}

In this section, it is shown that $H$ does not change under {\it some} Reidemeister moves that do not involve diagrams with more than $n$ crossings. It will be proven in latter sections that $H$ does not change under other moves. These other moves or {\it bad} moves are presented in Figure \ref{badmoves}. They involve basepoint or basepoints.

The bad moves $\Omega_4 '$, $\Omega_5 '$ and $\Omega_5 ''$ involve simple based diagrams.
The bad $\Omega_4 '$ move is the move in which, for some component, the secondary basepoints are changed.
The bad $\Omega_5 '$ move involves two couples of secondary basepoints and the bad $\Omega_5 ''$ move involves the primary and some secondary basepoints.

All other moves are {\it good} moves. For these moves the proof that $H$ is unchanged is similar to the case of Homfly for links in $\mathbb R^3$ (see \cite{LM}).

Using relation (HI) and induction IH($n-1$) one gets easily the following:

\begin{remark}[Crossing changes outside Reidemeister moves]\label{nice_remark}
Consider a Reidemeister move from a based diagram $\dot{D}_1$ to $\dot{D}_2$. If $H$ does not change under this move, than it does not change under the same move from $\dot{D}_1 '$ to $\dot{D}_2 '$ where $\dot{D}_1 '$ is obtained from $\dot{D}_1$ by switching some crossings not involved in the move, and $\dot{D}_2 '$ is obtained from $\dot{D}_2$ by switching the same crossings.
\end{remark}

\begin{lemma}[Crossing changes inside Reidemeister moves]\label{independence_reid}
Consider a Reidemeister move that is not an $\Omega_4$ move, from a based diagram $\dot{D}_1$ to $\dot{D}_2$.
Let $b_1$ and $b_2$ be two branches involved in the move. Suppose that they are not the lowest and uppermost branch in an $\Omega_3$ move.

If $H$ does not change under this move, then it does not change under another move obtained from the first one by switching the crossing(s) between $b_1$ and $b_2$.
\begin{proof}
In the case of an $\Omega_2$ move, an easy calculation using (HI) and (HII) shows that $H$ is unchanged, if one switches the two crossings that disappear under the move.

In the case of an $\Omega_3$ move, let $\dot{D}'_1$ (resp. $\dot{D}'_2$) be the diagram obtained from $\dot{D}_1$ (resp. $\dot{D}_2$) by switching two adjacent branches (for example the lowest and the middle ones).
Let $D''_1$ (resp. $D''_2$) be the diagram obtained from $\dot{D}_1$ (resp. $\dot{D}_2$) by smoothing the crossing that is switched to obtain $\dot{D'_1}$ (resp. $\dot{D'_2}$). By assumption $H(\dot{D_1})=H(\dot{D_2})$.
Now, $H(D''_1)=H(D''_2)$ because either the two diagrams are equal, or one can pass from one to the other by two $\Omega_2$ moves that do not increase the number of crossings beyond $n-1$ and one may use IH($n-1$).
From (HI) it follows that $H(\dot{D'_1})=H(\dot{D'_2})$.

The case of an $\Omega_5$ move is treated similarly to the case of an $\Omega_3$ move. One uses two $\Omega_4$ moves between diagrams with $n-1$ crossings and IH($n-1$) as well as (HI).
\end{proof}
\end{lemma}

\begin{lemma}[Invariance under good Reidemeister moves]
$H$ does not change under good Reidemeister moves that involve diagrams with at most $n$ crossings.
\begin{proof}
Because of IH($n-1$) it is sufficient to consider the case when at least one of the diagrams involved in the move has $n$ crossings.
If the basepoint of the diagram with $n$ crossings is lying on a 0-homologous component, it can be pushed out of the move by Lemma \ref{moving_basepoint} without changing $H$.

One may suppose that the diagram before the move has $n$ crossings and that it is descending using Remark \ref{nice_remark} and Lemma \ref{independence_reid}. It is easily checked that the diagram after the move is again descending except in a special case considered at the end of the proof. Now, from the definition of $H$ and IH($n-1$) it follows easily that $H$ is unchanged under the move.

The special case that has to be considered is the following: the basepoint is lying on a 0-homologous component and, after an $\Omega_2$ move, this component is not involved in any crossing. But again, in this case $H$ does not change by definition and IH($n-1$).

\end{proof}
\end{lemma}

\subsection{Basepoints of simple diagrams}
In this section it is shown that $H$ does not depend on the choice of antipodal basepoints for simple diagrams.

We say that an arc in a diagram $D$ {\it separates} a couple $P$, $Q$ of antipodal points lying on the boundary $S$ of the disk of $D$, if the endpoints of this arc are in different connected components of $S-(P\cup Q)$. Otherwise, we say that the arc {\it does not separate} the couple.

The following lemma is a reformulation of Lemma 1 of \cite{M}:
\begin{lemma}[Non separating arcs for components with at least three arcs]\label{old_lemma}
Let $D$ be a simple diagram. Let $b$ be a 1-homologous component. Let $P$, $Q$ be a couple of antipodal points on the boundary circle of $D$, $P$ and $Q$ not in $b$. Suppose that $b$ has at least three arcs.

Then at least two of the arcs of $b$ do not separate $P$, $Q$.
\end{lemma}

\begin{lemma}[Non separating arcs for components with at least five arcs]\label{4points}
Let $D$ be a simple diagram and $b$ a 1-homologous component of $D$.
Suppose that $b$ has at least five arcs.
Consider two couples of antipodal endpoints of some arcs of $b$, such that moving from one couple to another in the counterclockwise direction, no other endpoints of $b$ are encountered.

Then there is at least one arc with no endpoint in these two couples, which does not separate any of the two couples.
\begin{proof}

In Figure \ref{4cases} four cases are presented. $P^1$ and $Q^1$ is a couple of antipodal endpoints. $P^2$ and $Q^2$ is another such couple.
Suppose that there are no endpoints of $b$ between $P^1$ and $P^2$ (thus no endpoints between $Q^1$ and $Q^2$).

\begin{figure}[ht]
\scalebox{1}{\includegraphics{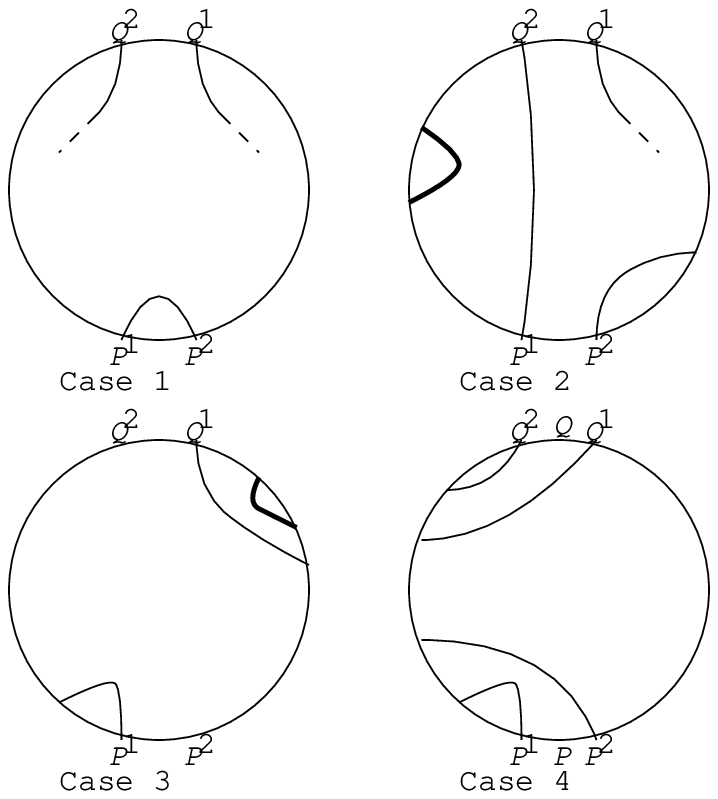}}
\caption{}
\label{4cases}
\end{figure}

{\it Case 1}: $P^1$ and $P^2$ are endpoints of the same arc. The $\Omega_4$ move, which deletes this arc, gives rise to a simple diagram in which $b$ has at least three arcs.
This move changes only the three arcs shown in Figure \ref{4cases}, Case 1.
By Lemma \ref{old_lemma}, at least two arcs of $b$ do not separate the antipodal endpoints $P^1$ and $Q^1$ (and also $P^2$ and $Q^2$). One of these may coincide with the arc obtained from the three arcs appearing in Case 1. But there is an extra arc not separating the couples of endpoints.

{\it Case 2}: Suppose that $P^1$ and $Q^2$ are joined by an arc. Then the required arc exists obviously (the one drawn in thick).

{\it Case 3}: Suppose that the arcs starting at $P^1$ and $Q^1$ do not separate the couple $P^2$ and $Q^2$. Then one finds again a suitable arc (drawn in thick pencil).

{\it Case 4}: In the only possibility which still has to be considered, there has to be some arcs with endpoints $P^1$ or $Q^1$ that separate the couple $P^2$ and $Q^2$ or vice versa.
Suppose for instance, that the arc which has $P^2$ as endpoint separates the couple $P^1$ and $Q^1$. Then the arc with endpoint $P^1$ cannot separate the couple $P^2$ and $Q^2$ (the diagram is simple). Since we may assume that this is not a situation of Case 3, the arc with endpoint $Q^1$ has to separate the couple $P^2$ and $Q^2$. And this finally means that the arc with endpoint $Q^2$ does not separate the couple $P^1$ and $Q^1$.

Denote by $P$ and $Q$ a couple of antipodal points on the boundary circle of the diagram, that are not endpoints, such that $P$ is between $P^1$ and $P^2$ and $Q$ is between $Q^1$ and $Q^2$.

Now we have two arcs which do not separate the couple $P$ and $Q$ (namely, arcs with endpoints $P^1$ and $Q^2$) whereas two arcs separate this couple (the arcs with endpoints $P^2$ and $Q^1$). Note that the two arcs that do not separate $P$, $Q$ are on the left side of the diagram.

To find a suitable arc in this case, let us prove first that the number of arcs of $b$ not separating $P$ and $Q$ on the left side is the same as on the right side. This can be checked by induction on the number of arcs of $b$: if $b$ has a unique arc this number is 0 on both sides. Otherwise, if one uses an $\Omega_4$ move to remove an arc of $b$ that does not separate $P$, $Q$ (which exists by Lemma \ref{old_lemma}), the number of arcs not separating this couple on the left and on the right is decreased by one, or remains constant.

Thus there are at least two arcs of $b$ on the right side side of the diagram that do not separate $P$, $Q$ and thus not separating $P_1$, $Q_1$ and $P_2$, $Q_2$ as well.

\end{proof}
\end{lemma}

A {\it subarc} is a compact connected submanifold of an arc.

Let $a$ be a subarc going from a crossing to itself, such that it has no other self-crossings. Let $P$ be a point on the arc of which $a$ is a subarc, just outside $a$.
Then $a$ is called a {\it 1-gon} if, in the net, it does not separate $P$ and the line at infinity. 

Let ($a_1$, $a_2$) be a couple of subarcs both going from one crossing to another, having no extra crossings between them, and none of them having self-crossings. Let $P_1$ be a point on an arc of which $a_1$ or $a_2$ is a subarc, just outside $a_1$ and $a_2$, and close to the first crossing between $a_1$ and $a_2$. Let $P_2$ be a point with the same properties as $P_1$ except that it is close to the second crossing between $a_1$ and $a_2$.
Then ($a_1$, $a_2$) is called a {\it 2-gon} if, in the net, it does not separate $P_1$ and the line at infinity and it does not separate $P_2$ and the line at infinity.

\begin{figure}[ht]
\scalebox{1}{\includegraphics{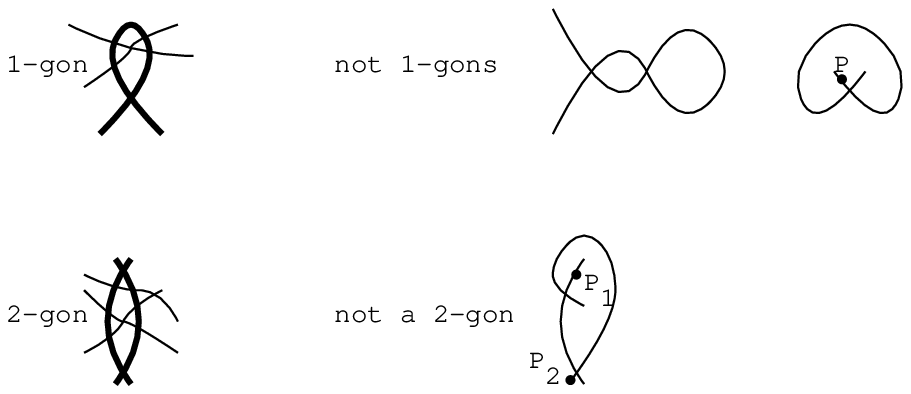}}
\caption{}
\label{loop_and_fish}
\end{figure}

1-gons and 2-gons are presented in Figure \ref{loop_and_fish}.

By definition, a subarc $a$ is {\it below} a subarc $b$, if at each crossing involving $a$ and $b$ the branch in $a$ is under the branch in $b$.

\begin{lemma}[moving subarcs]\label{moving_arcs}
Let $D$ be a (possibly based) diagram with at most $n$ crossings.
Let $B$ be a disk in $D$ such that there is no basepoint in $B$, and such that $B$ does not intersect the boundary circle of the diagram.
Suppose that $B\cap D$ is the union of a subarc $a$ strictly included in $\partial B$ and some other subarcs $a_1, a_2, ..., a_l$ properly embedded in $B$.
Suppose that $a_i$ is below $a_j$ for any $i<j$ and that there is some $k$ such that $a$ is above $a_i$ for $i\le k$ and $a$ is below $a_i$ for $i>k$.

Suppose that each $a_i$ crosses $a$ in one point and that no pair of $a_i$-s has more than one common crossing. Let $b$ be the closure of $\partial B-a$ and let $D'$ be the result of substituting $b$ for $a$ in $D$, where $b$ is above $a_i$ if $a$ is above $a_i$ and below otherwise. 

Then $H(D)=H(D')$.
\begin{proof}
It is easily seen that it is possible to transform $D$ into $D'$ with a series of good $\Omega_3$ moves.
\end{proof}
\end{lemma}

\begin{lemma}[removing 2-gons]\label{removing_2gons}
Suppose that $B$ and $D$ are as in Lemma \ref{moving_arcs}, except that one subarc $a_i$ crosses $a$ in two points, and each $a_j$, $j\ne i$, crosses $a_i$ in one point.
If, as before, $D'$ is the result of substituting $b$ for $a$ in $D$, then $H(D)=H(D')$.
\begin{proof}
Applying Lemma \ref{moving_arcs}, move all crossings in the interior of the 2-gon formed by $a$ and $a_i$ out of this 2-gon. Transform the resulting diagram into $D'$ by a series of good $\Omega_3$ moves.
\end{proof}
\end{lemma}

\begin{lemma}[removing 1-gons]\label{removing_1gons}
Let $D$ be a (possibly based) diagram with at most $n$ crossings.
Let $a$ be a 1-gon, bounding a disk $B$. Suppose that there are no 1-gons and no basepoint in $B$ (except for, possibly, the self-crossing of $a$ which can be a basepoint).
Suppose that $B\cap D$ is the union of $a$ and some other subarcs $a_1, a_2, ..., a_l$ properly embedded in $B$.
Suppose that $a_i$ is below $a_j$ for any $i<j$ and that there is some $k$ such that $a$ is above $a_i$ for $i\le k$ and $a$ is below $a_i$ for $i>k$.

Let $\epsilon$ be the sign of the self-crossing of $a$ and $D'$ the diagram obtained from $D$ by removing the 1-gon $a$.
Then $H(D)=(xv^{-1})^{\epsilon} H(D')$.
\begin{proof}
First, remove all 2-gons that are inside $B$ using Lemma \ref{removing_2gons}, starting with the most nested ones.
Then decrease the number of crossings inside $B$ using Lemma \ref{removing_2gons} for couples of subarcs in which one subarc is part of $a$. In this way, the number of crossings inside $B$ is reduced to 0.
Finally use relation (HII) which holds for diagrams with at most $n$ crossings (Proposition \ref{relation}).
\end{proof}
\end{lemma}

\begin{lemma}[special case of invariance under $\Omega_5'$ move]\label{five_prime}
Suppose that an $\Omega_5'$ move is applied on a simple based diagram $\dot{D}$ with $n$ crossings. Suppose that each of the two 1-homologous components involved in this move consists of a unique arc. Let $\dot{D}'$ be the diagram after the move.

Then $H(\dot{D})=H(\dot{D}')$.
\begin{proof}
By Remark \ref{nice_remark} and Lemma \ref{independence_reid} one may suppose that $\dot{D}$ is descending.
If the two components involved in the move have at least three common crossings, then, using Lemma \ref{removing_2gons}, remove two crossings that are not involved in the move without changing $H$ and, using IH($n-1$), get $H(\dot{D})=H(\dot{D}')$. Otherwise, if the two components have a unique common crossing, then $\dot{D}'$ is also descending, and both $H(\dot{D})$ and $H(\dot{D}')$ are equal to the same $H(\dot{D}_{k,l})$ for some $k$ and $l$.
\end{proof}
\end{lemma}

In the figures that follow, a sign of equality between two diagrams means that $H$ is the same for both.

\begin{proposition}[special case of independence on basepoints for simple diagrams]\label{simple_unchanged_special}
Let $\dot{D}_1$ and $\dot{D}_2$ be two simple based diagrams with $n$ crossings. Suppose that all 1-homologous components in $\dot{D}_1$ and $\dot{D}_2$ consist of unique arcs.
Suppose that $\dot{D}_1$ and $\dot{D}_2$ differ only by the position of the basepoints.

Then $H(\dot{D_1})=H(\dot{D_2})$.
\begin{proof}
It is sufficient to prove that $H(\dot{D}_1)=H(\dot{D}_2)$ if the basepoints of $\dot{D}_2$ are next to the basepoints of $\dot{D}_1$ in the counterclockwise direction.

By Lemma \ref{crossing_independence}, one may suppose that $\dot{D}_1$ is descending.

If there is a couple of components which have at least three crossings between them, then $H(\dot{D}_1)=H(\dot{D}_2)$: indeed, one may reduce the number of crossings in $\dot{D}_1$ and $\dot{D}_2$ using Lemma \ref{removing_2gons} (because  $\dot{D}_1$ is descending) then use IH($n-1$).
Suppose now that any two components have a unique common crossing.

In changing the position of a couple of antipodal basepoints from $\dot{D}_1$ to $\dot{D}_2$, the following two cases can occur:

$\bullet$ {\it Case 1}: Traveling from the basepoint of $\dot{D}_1$ that is initial, in the counterclockwise direction, the first secondary basepoint encountered is also initial (see Figure \ref{supersimple}).

One may suppose that $\dot{D}_1$ is in fact some $\dot{D}_{k,l}$. Indeed, if it is not, use several times Lemma \ref{moving_arcs} and $\Omega_5'$ move on components with secondary basepoints (Lemma \ref{five_prime}). Then, using Lemma \ref{crossing_independence}, one may suppose that $\dot{D}_1$ is equal to $\alpha(\dot{D}_{k,l})$, i.e. it is almost standard and the sign at every crossing is equal to +1.

By $(H_3)$, $H(\dot{D_1})=z^{k+l}$.
As it is shown in Figure \ref{supersimple}, $H(\dot{D}_1)=H(\dot{D}_2)$.

\begin{figure}[ht]
\scalebox{1}{\includegraphics{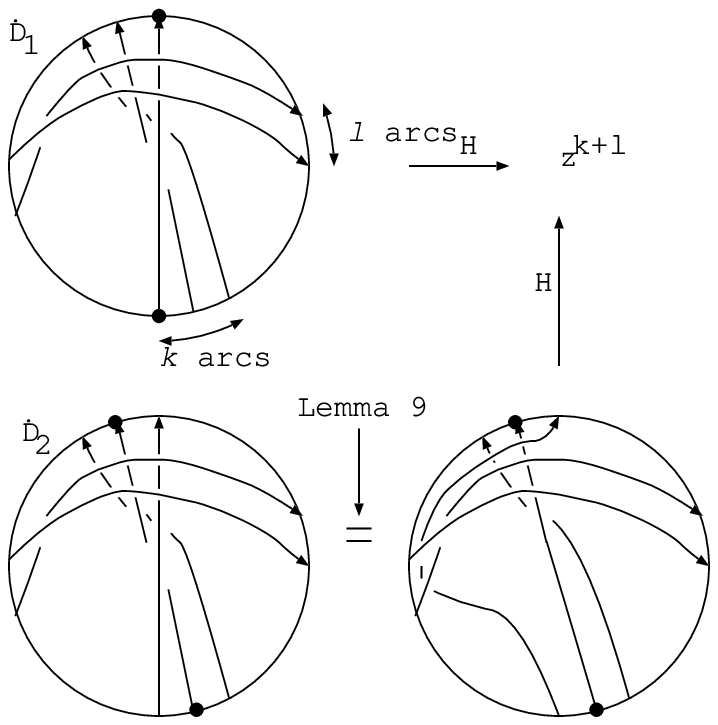}}
\caption{}
\label{supersimple}
\end{figure}

In order to prove Case 2, we need to have invariance of $H$ under some special $\Omega_5 ''$ move.
Suppose that $\dot{D}$ is a simple based diagram with $n$ crossings such that each of its 1-homologous components consists of a unique arc.
Suppose that one applies on $\dot{D}$ an $\Omega_5$ move involving primary and secondary basepoints (bad $\Omega_5''$ move), and the vanishing triangle has two vertices that are both initial basepoints, or both final endpoints. Suppose also that the two components involved in the move have a unique common crossing (it is the crossing that appears in the move).

Then $H$ does not change under such move. Indeed, using Remark \ref{nice_remark} and Lemma \ref{independence_reid}, suppose that $\dot{D}$ is descending. Before the move $H(\dot{D})=H(\dot{D}_{k,l})$ for some $k$ and $l$. After the move, change the basepoints in clockwise direction: as it was just proven above, $H$ does not change. The resulting based diagram is again descending and $H$ of this diagram is again equal to $H(\dot{D}_{k,l})$.

$\bullet$ {\it Case 2}: Traveling from the basepoint of $\dot{D}_1$ that is initial, in the counterclockwise direction, the first secondary basepoint encountered is final (see Figure \ref{supersimplebis}).

Suppose that $\dot{D}_1$ is as on the top, left in Figure \ref{supersimplebis}. One may always reduce Case 2 to the situation presented on this figure, using Lemma \ref{crossing_independence}, Lemma \ref{moving_arcs} and Lemma \ref{five_prime}.

\begin{figure}[ht]
\scalebox{1}{\includegraphics{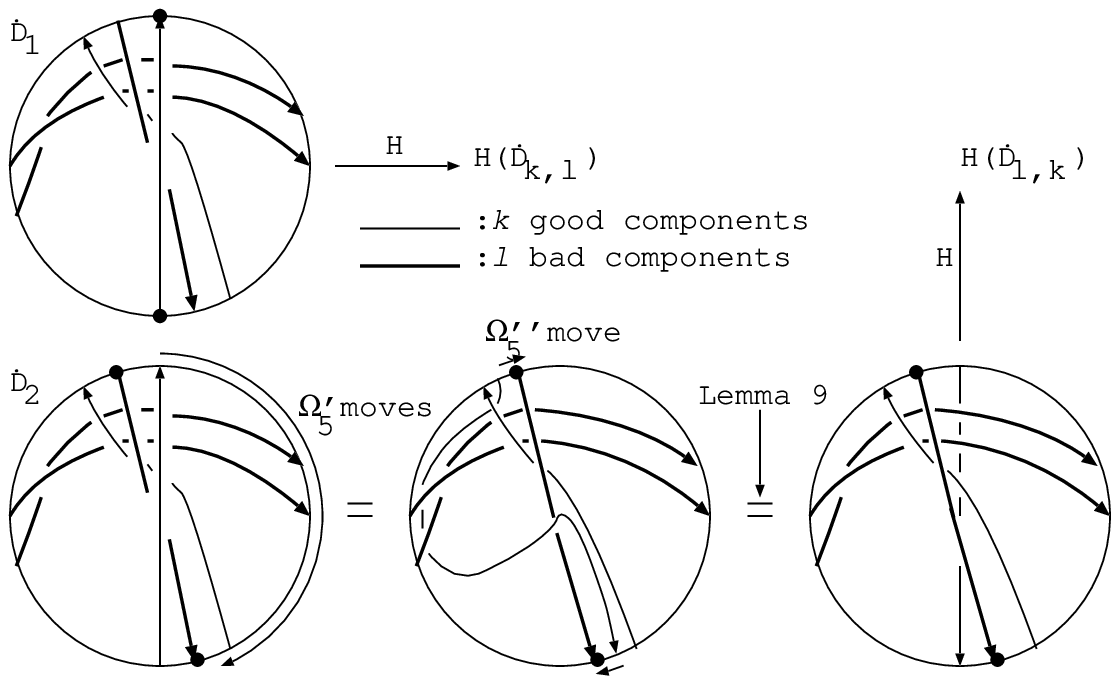}}
\caption{}
\label{supersimplebis}
\end{figure}

In this case, $H(\dot{D}_1)=H(\dot{D}_{k,l})$, whereas $H(\dot{D}_2)=H(\dot{D}_{l,k})$.
But, as it is shown in Figure \ref{kllk}, $H(\dot{D}_{k,l})=H(\dot{D}_{l,k})$. In this figure bad $\Omega_5''$ moves are applied for which the vanishing triangles have two vertices that are both initial or both final basepoints. Also, {\it Case 1} is used to change the position of primary basepoints without changing $H$.

\begin{figure}[ht]
\scalebox{1}{\includegraphics{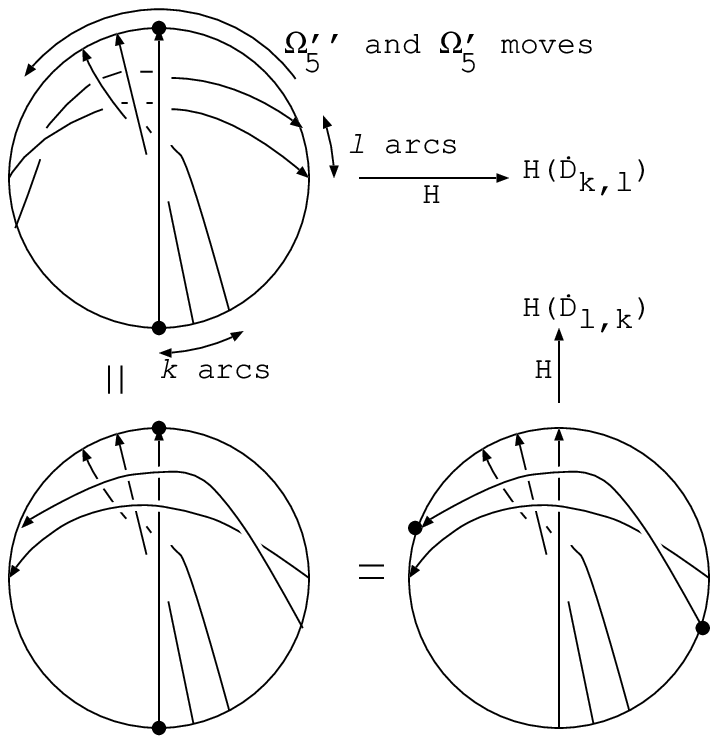}}
\caption{}
\label{kllk}
\end{figure}

\end{proof}
\end{proposition}

\begin{lemma}[special case of invariance under $\Omega_5''$ move]\label{five_bis}
Suppose that $\dot{D}$ is a simple based diagram with $n$ crossings. Suppose that all 1-homologous components of $\dot{D}$ consist of unique arcs.
Let $\dot{D}'$ be the diagram obtained from $\dot{D}$ by the application of an $\Omega_5''$ move.

Then $H(\dot{D})=H(\dot{D}')$.
\begin{proof}

Using Remark \ref{nice_remark} and Lemma \ref{independence_reid}, assume that $\dot{D}$ is descending.

First, suppose that the two components involved in the move have at least three common crossings. Then $H(\dot{D})=H(\dot{D}')$, because one can remove two crossings by Lemma \ref{removing_2gons} and use IH($n-1$). We may therefore assume that the two components involved in the move have a unique common crossing.

If, in applying the $\Omega_5''$, the vanishing triangle has two vertices that are both initial basepoints, or both final endpoints then $H(\dot{D})=H(\dot{D}')$ as it was already seen in the proof of Proposition \ref{simple_unchanged_special}.

Suppose now that the vanishing triangle has one vertex that is an initial basepoint and the other one that is a final basepoint. As $\dot{D}$ is descending, $H(\dot{D})=H(\dot{D}_{k,l})$ for some $k$ and $l$.
By moving the basepoints of $\dot{D}'$ in the clockwise direction, one gets a descending diagram with $l$ good and $k$ bad components. Thus $H(\dot{D}')=H(\dot{D}_{l,k})$. But it was seen in the proof of Proposition \ref{simple_unchanged_special} that $H(\dot{D}_{k,l})=H(\dot{D}_{l,k})$.
\end{proof}
\end{lemma}

\begin{proposition}[independence on basepoints for simple diagrams]\label{simple_unchanged}
Let $\dot{D}_1$ and $\dot{D}_2$ be two simple based diagrams with $n$ crossings. Suppose that $\dot{D}_1$ and $\dot{D}_2$ differ only by the position of the basepoints.
Then $H(\dot{D_1})=H(\dot{D_2})$.

\begin{proof}
It is sufficient to prove that $H(\dot{D}_1)=H(\dot{D}_2)$ if the basepoints of $\dot{D}_2$ are next to the basepoints of $\dot{D}_1$ in the counterclockwise direction.
Let $(P_1,Q_1)$ (resp. $(P_2,Q_2)$) be the couple of basepoints of $\dot{D}_1$ (resp. $\dot{D}_2$). Suppose that $(P_2,Q_2)$ is next to $(P_1,Q_1)$ in the counterclockwise direction.

By Lemma \ref{crossing_independence}, we may assume that $\dot{D}_1$ is descending.

Let $a$ be the component to which $(P_1,Q_1)$ belongs. Then $(P_2,Q_2)$ may belong to the same component $a$ or to a different one, say $b$.

$\bullet$ {\it Case 1}: Suppose that $(P_2,Q_2)$ belongs to $a$. 

Note that the secondary basepoints of $\dot{D}_2$ coincide with the secondary basepoints of $\dot{D}_1$. Thus, the ordering of components of $\dot{D}_1$ arising from the primary basepoints is the same for $\dot{D}_2$.
Consider, for this ordering, the last component that has more than one arc, if there is such component. In Figure \ref{atoa}, it is the component that is dashed.

\begin{figure}[ht]
\scalebox{1}{\includegraphics{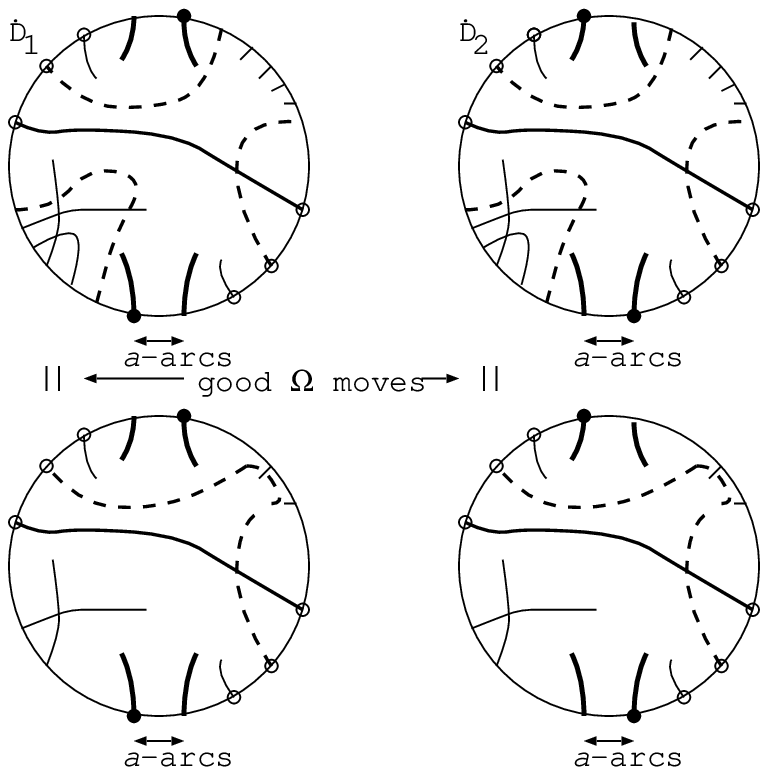}}
\caption{}
\label{atoa}
\end{figure}

One can reduce the number of arcs of this component without changing $H$ for both $\dot{D}_1$ and $\dot{D}_2$. By Lemma \ref{old_lemma}, it has an arc, say $c$, that does not separate its secondary couple of endpoints.
Consider an arc that is most nested in it. This arc cannot have an endpoint that is a secondary basepoint.
Using Lemma \ref{moving_arcs}, Lemma \ref{removing_2gons} and good $\Omega_5$ and $\Omega_4$ moves (good in $\dot{D}_1$ {\it and} in $\dot{D}_2$), this arc is removed without changing $H$. The diagram obtained in this way from $\dot{D}_1$ is still descending. In the same manner, all arcs nested in the arc $c$ are removed. Finally $c$ is removed. 

In this way, without changing $H$, $\dot{D}_1$ is transformed into a descending based diagram (still denoted by $\dot{D}_1$) in which all components, except $a$, have unique arcs. $\dot{D}_2$ is transformed similarly.
Now using Lemma \ref{4points} the number of arcs of $a$ is reduced to three, by eliminating the arcs that do not separate the couples ($P_1$,$Q_1$) and ($P_2$,$Q_2$), and that have no endpoints coinciding with $P_1$, $Q_1$, $P_2$ or $Q_2$. During the elimination of such arcs, $H$ is unchanged for $\dot{D}_1$ and $\dot{D}_2$, and $\dot{D_1}$ stays descending.
One arrives at the situation presented in Figure \ref{atoa2}. In this figure the arcs of $a$ (the thickest ones) are marked with {\it high} if they are above everything else, or {\it low} if they are below everything else.

\begin{figure}[ht]
\scalebox{1}{\includegraphics{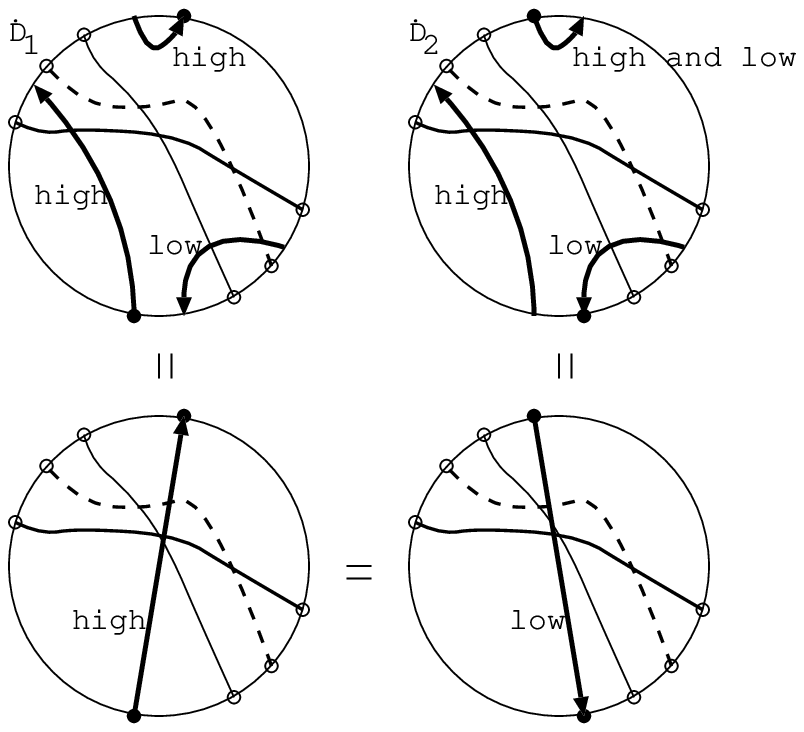}}
\caption{}
\label{atoa2}
\end{figure}

At the bottom of the figure, the equality holds because one may turn around the component $a$ by moving some subarcs (Lemma \ref{moving_arcs}) and by application of several bad $\Omega_5''$ moves in the situation where all components have unique arcs. It follows from Lemma \ref{five_bis} that $H$ does not change.

$\bullet$ {\it Case 2}: Suppose that $(P_2,Q_2)$ belongs to a component $b$, different from $a$.

A first possibility is that $a$ has a unique arc. In that case, in the same way as in Case 1, one can reduce without changing $H$ the number of arcs and get diagrams for which all components have unique arcs. Then, it follows from Proposition \ref{simple_unchanged_special} that $H$ is unchanged.

Otherwise $a$ has several arcs. The method of {\it Case 1} can be repeated for components that have secondary basepoints after the secondary basepoints of $a$ in $\dot{D}_2$. For each of these components the number of arcs is reduced to one. Also the number of arcs of $a$ is reduced to three with the help of Lemma \ref{4points}.
Now, one has to compare the diagrams at the top of Figure \ref{a_to_b}.
In this figure, the arcs of $a$ (the thickest ones) are again marked with {\it high} or {\it low}. Note that in the diagram on the right bottom of the figure, the unique arc of $a$ is divided into two parts: one {\it high} and the other one {\it low}.

\begin{figure}[ht]
\scalebox{1}{\includegraphics{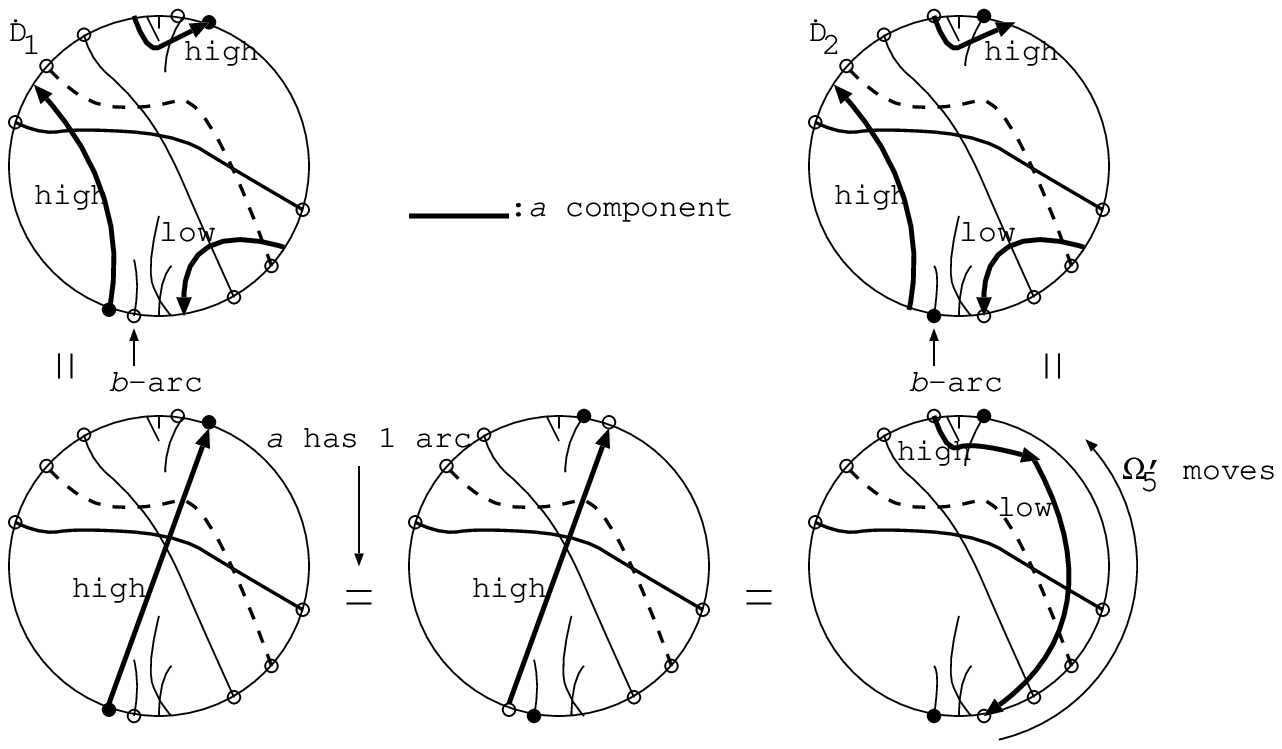}}
\caption{}
\label{a_to_b}
\end{figure}

In the diagram on the right bottom of the Figure \ref{a_to_b}, one applies several bad $\Omega_5'$ moves involving $a$, which has a unique arc, and another component with unique arc. By Lemma \ref{five_prime}, $H$ does not change under such move. In general, one applies also good $\Omega_5$ moves involving $a$ and some components that have secondary basepoints before the secondary basepoints of $a$ in $\dot{D}_2$. One may also have to eliminate some arcs of these components using Lemma \ref{moving_arcs} and good $\Omega_5$ and $\Omega_4$ moves.

Finally, $H$ is the same on the left and right side of Figure \ref{a_to_b}.
\end{proof}
\end{proposition}

\subsection{Bad moves $\Omega_4', \Omega_5'$ and $\Omega_5''$}

\begin{proposition}[invariance under $\Omega_4', \Omega_5'$ and $\Omega_5''$ moves]
$H$ does not change under bad moves $\Omega_4', \Omega_5'$ and $\Omega_5''$, involving diagrams with $n$ crossings.
\begin{proof}
Recall that a bad $\Omega_4'$ move is an $\Omega_4$ move under which the secondary basepoints of some component are changed.
An $\Omega_4'$ move is shown at the top of Figure \ref{bad4}. From this figure it is clear that $H$ does not change under the move, because it can be transformed to a good $\Omega_4$ move, by changing the position of the primary basepoints.

\begin{figure}[ht]
\scalebox{1}{\includegraphics{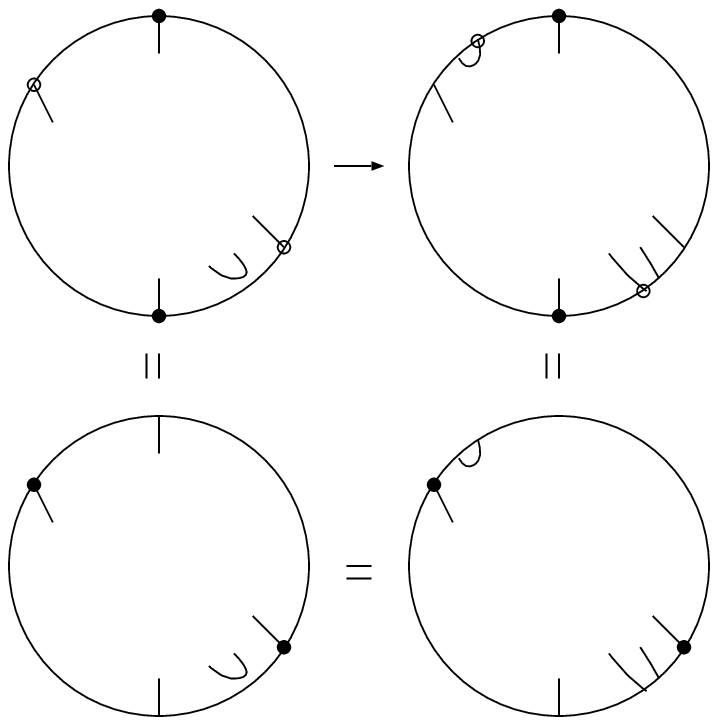}}
\caption{}
\label{bad4}
\end{figure}

A bad $\Omega_5'$ move involving two components which have unique arcs leaves $H$ unchanged according to Lemma \ref{five_prime}.
Now, if a component involved in $\Omega_5'$ move has at least three arcs then, by changing the position of the basepoints, the move can be transformed into a good $\Omega_5$ move. Thus any $\Omega_5'$ move leaves $H$ unchanged.

If all components have unique arcs the bad $\Omega_5''$ move leaves $H$ unchanged according to Lemma \ref{five_bis}.
Now if there is a component with at least three arcs, one can transform the $\Omega_5''$ move into an $\Omega_5$ or $\Omega_5'$ move, by changing the position of the basepoints. It follows that any $\Omega_5''$ move leaves $H$ unchanged.
\end{proof}
\end{proposition}

\subsection{Bad moves $\Omega_2'$ and $\Omega_3'$}
In this section, the invariance of $H$ under bad moves $\Omega_2'$ and $\Omega_3'$ is established as a consequence of Proposition \ref{shortening} below.

\begin{lemma}[removing 1-gons, stronger version]\label{removing_1gons_strong}
Let $D$ be a diagram with at most $n$ crossings.
Let $a$ be a 1-gon, bounding a disk $B$. Suppose that there is no basepoint in $B$ (except, possibly, for the self-crossing of $a$ which can be a basepoint).
Suppose that $B\cap D$ is the union of $a$ and some other subarcs $a_1, a_2, ..., a_l$ properly immersed in $B$.
Suppose that there is some $k$ such that $a$ is above $a_i$ for $i\le k$ and $a$ is below $a_i$ for $i>k$. Suppose also that if $i\le k$ and $j>k$ then $a_i$ is below $a_j$.

Let $\epsilon$ be the sign of the self-crossing of $a$ and $D'$ the diagram obtained from $D$ by removing the 1-gon $a$.
Then $H(D)=(xv^{-1})^{\epsilon} H(D')$.
\begin{proof}
The proof is done by induction on the number of crossings in the interior of $B$, say $m$. If $m=0$ then one may apply Lemma \ref{removing_1gons}.

If there are no 1-gons in the interior of $B$, use relation (HI) on both $D$ and $D'$ to order the subarcs $a_1, ..., a_l$ so that $a_i$ is below $a_j$ if $i<j$. Using (HI) gives rise to smoothings for which one applies induction on $m$.
Now for the diagrams in which $a_1, ..., a_l$ are ordered, one applies Lemma \ref{removing_1gons}.

If there are some 1-gons in the interior of $a$, consider one of them that is most nested (i.e. there are no 1-gons in its interior). It can be eliminated as in the preceding paragraph and one applies induction on $m$.
\end{proof}
\end{lemma}

\begin{lemma}[removing triangles]\label{triangle}
Suppose that in a diagram $D$ with at most $n-1$ crossings, one has the situation presented on the left of Figure \ref{triangle1}. At the bottom of $D$ two subarcs $a$ and $b$ (the thick ones), together with a part of the boundary circle of $D$, form a triangle. Inside the triangle, there are properly embedded subarcs $a_1, a_2, ..., a_l$, all intersecting $a$ and $b$ in a unique crossing. Suppose that no pair of $a_i$-s has more than one crossing.
Suppose that $a_i$ is below $a_j$ if $i<j$ and that there are $k$ and $k'$ such that $a$ is above $a_i$ for $i\le k$, $a$ is below $a_i$ for $i>k$, $b$ is above $a_i$ for $i\le k'$, $b$ is below $a_i$ for $i>k'$. Suppose that $a$ is below $b$ if $k\le k'$ and $a$ is above $b$ otherwise.

Then $H$ is the same for diagrams on the left and on the right of Figure \ref{triangle1}, where $a$, $b$ and the $a_i$-s are above or below each other on the right in the same way as they are on the left.

\begin{figure}[ht]
\scalebox{1}{\includegraphics{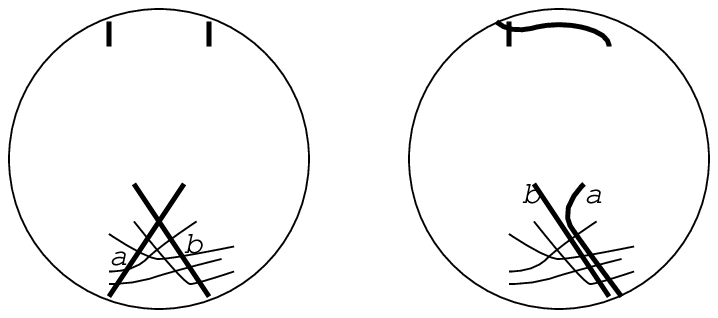}}
\caption{}
\label{triangle1}
\end{figure}
\begin{proof}
The proof is similar to the proof of Lemma \ref{moving_arcs}. The diagram on the left of Figure \ref{triangle1} is transformed into the diagram on the right of this figure, by application of several $\Omega_3$ moves and one $\Omega_5$ move. This transformation is sketched in Figure \ref{triangle2}.

\begin{figure}[ht]
\scalebox{1}{\includegraphics{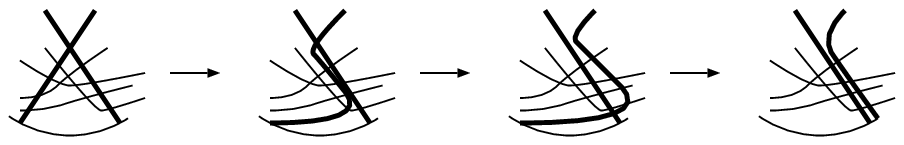}}
\caption{}
\label{triangle2}
\end{figure}
\end{proof}
\end{lemma}

\begin{proposition}[shortening of diagrams]\label{shortening}
Let $D$ be a diagram with $k$ crossings. A part of $D$ is shown on the left of Figure \ref{lemma_desc}.
Suppose that the arc distance from $P$ to $Q$ is even, and suppose that $D$ is descending from $P$ to $Q$.

The part of $D$ from $P$ to $Q$ is dashed.
$D'$ is obtained from $D$ by removing this dashed part, and joining $P$ and $Q$ with a segment (a part of $D'$ is shown on the right of Figure \ref{lemma_desc}). Let $w$ be the sum of signs at all crossings for which both branches are in the dashed part.

If $k\le n-1$ then $H(D)=(xv^{-1})^w H(D')$.

\begin{figure}[ht]
\scalebox{1}{\includegraphics{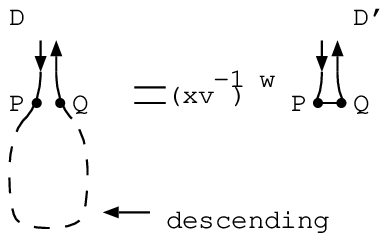}}
\caption{}
\label{lemma_desc}
\end{figure}

\begin{proof}
The proof is done by induction on $k$, the number of crossings in $D$.
If $k=0$ then $H(D)=H(D')$ by definition of $H$.

Suppose the proposition is true for $k<l$ where $l\le n-1$. Suppose now that $k=l$.
We will use induction on the arc distance from $P$ to $Q$. The case when this distance is equal to $0$ (the dashed part consists of a unique subarc) will be considered at the end of the proof.

Suppose now that the arc distance from $P$ to $Q$ is at least two.
This means that there are several arcs that are dashed or partially dashed.
The arc that one encounters first when traveling in the net from $P$ and crossing the line at infinity once is below all other arcs.
It may have 1-gons if there is a crossing with both branches in this lowest arc, or not.

$\bullet$ {\it Case 1}: there are 1-gons in the lowest arc (see Figure \ref{lemma_multiarc}).

\begin{figure}[ht]
\scalebox{1}{\includegraphics{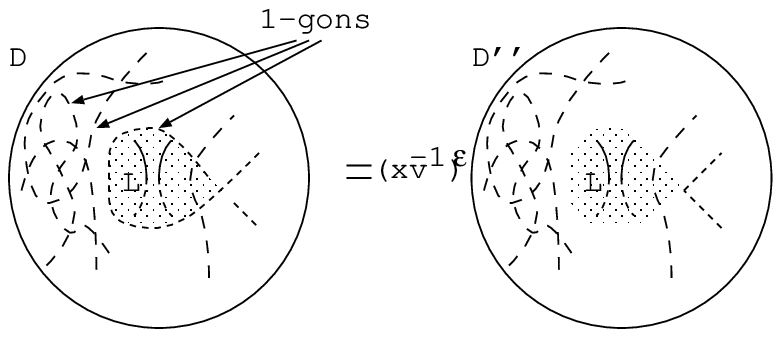}}
\caption{}
\label{lemma_multiarc}
\end{figure}

Consider $D''$ obtained from $D$ by eliminating some chosen 1-gon in the lowest arc. Let $\epsilon$ be the sign at the self-crossing of this 1-gon.
Notice that $H(D'')=(xv^{-1})^{w-\epsilon } H(D')$ by induction on $k$ ($D''$ has strictly less crossings than $D$).

To see that one has equality between $H(D)$ and $(xv^{-1})^\epsilon H(D'')$, one considers the part $L$ of $D$ that is inside the 1-gon (and the corresponding part of $D''$). In $L$, there can be some subarcs below the 1-gon (such subarcs can only be parts of the lowest arc). Otherwise, all other subarcs are above the 1-gon, and they are above the subarcs that are below the 1-gon.
Lemma \ref{removing_1gons_strong} can be applied, so $H(D)=(xv^{-1})^\epsilon H(D'')$.
Thus $H(D)=(xv^{-1})^{w} H(D')$.

$\bullet$ {\it Case 2}: there are no 1-gons in the lowest arc.

The lowest arc divides $D$ in two parts. In one of these parts are the two points $P$ and $Q$. Consider the other part, the {\it good} part.

There may be 1-gons that are entirely inside the good part.
If there are such 1-gons that are dashed, one gets $H(D)=(xv^{-1})^w H(D')$ using Lemma \ref{removing_1gons_strong} as in Case 1 (here it is essential that $P$ and $Q$ are not inside such 1-gons).

If inside the good part there are only 1-gons that are not dashed, then, using for $D$ and $D'$ relation (HI) on crossings that have both branches not in the dashed part, one gets again to a situation where an application of Lemma \ref{removing_1gons_strong} is possible, whereas for the smoothings that appear in (HI), one uses induction on $k$.

Suppose now, that there are no 1-gons in the good part of the lowest arc.

Inside this part consider a most nested dashed arc. Two situations can occur, a most nested arc can have a good part which either contains a couple of antipodal points on the boundary circle of the diagram (more difficult case), or does not contain such a couple (simpler case). These two cases are shown in Figure \ref{lemma_arcs}.

\begin{figure}[ht]
\scalebox{1}{\includegraphics{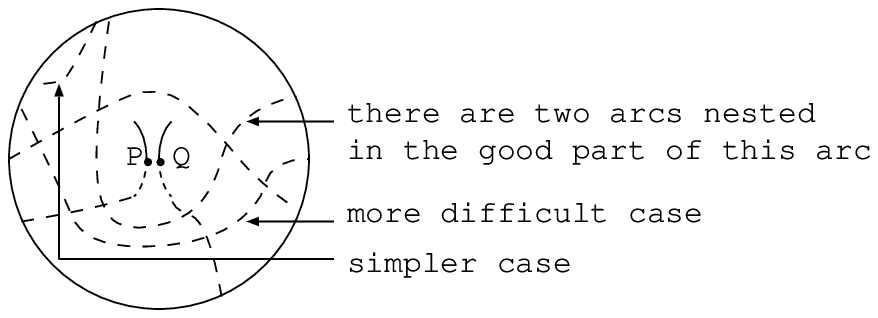}}
\caption{}
\label{lemma_arcs}
\end{figure}

In the simpler case $H(D)$ and $H(D')$ are computed using relation (HI), by ordering one above the other, the non dashed subarcs inside the good part of the most nested dashed arc. Again, for smoothings one has equality for $H$ (up to $(xv^{-1})^w$) by induction on $k$.
For the diagrams where the subarcs are ordered, one gets equality for $H$ (up to $(xv^{-1})^w$) using either Lemma \ref{removing_2gons} and induction on $k$; or Lemma \ref{moving_arcs}, $\Omega_5$ and $\Omega_4$ moves to make the dashed arc disappear and induction on the arc distance from $P$ to $Q$.

In the more difficult case the situation is presented in Figure \ref{difficult}.

\begin{figure}[ht]
\scalebox{1}{\includegraphics{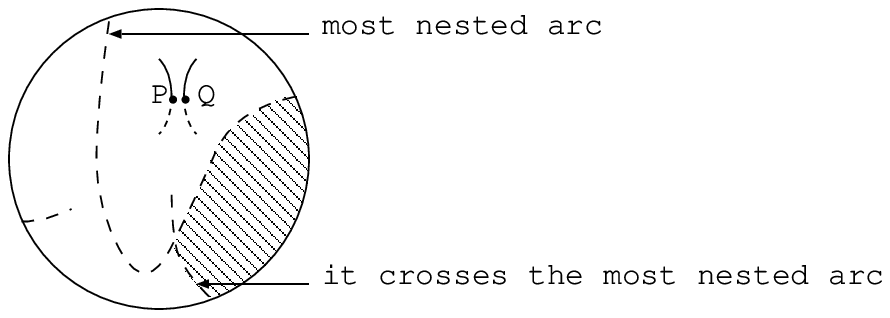}}
\caption{}
\label{difficult}
\end{figure}

Consider an arc which has an endpoint antipodal to an endpoint of the most nested dashed arc. Such arc has to intersect this most nested arc: otherwise it would be nested in this most nested arc.

Consider the dashed triangle in Figure \ref{difficult}. Order all subarcs inside the triangle with (HI) relation (for the smoothings use induction on $k$ to get equality for $H$ up to $(xv^{-1})^{w})$.
If there are 2-gons inside the triangle (including 2-gons formed by a part of a side of the triangle and another subarc), eliminate them without changing $H$ (Lemma \ref{removing_2gons}) and decrease the number of crossings. Then use induction on $k$.

Now, suppose that there are no 2-gons inside the triangle. First, eliminate all subarcs inside the triangle with both endpoints in the boundary circle of the diagram (i.e. all arcs inside the triangle) using Lemma \ref{moving_arcs}, $\Omega_5$ and $\Omega_4$ moves.

If the side of the triangle, which is part of the boundary circle of the diagram, contains no endpoints of arcs, except for the two vertices of the triangle, then, by Lemma \ref{triangle}, the triangle can be removed without changing $H$ and one gets to the simpler case.

If the side of the triangle, which is part of the boundary circle of the diagram, contains endpoints of arcs that are not vertices of the triangle, consider a smaller triangle inside the original one, which has one side in the boundary circle of the diagram and for which the assumptions of Lemma \ref{triangle} are satisfied. It is easily seen that such a triangle can always be found.
The smaller triangle is eliminated without changing $H$.
After the elimination the number of crossings inside the original triangle {\it or} the number of endpoints of arcs in the original triangle is decreased.
Repeating this procedure several times, one gets to a situation where there are no endpoints of arcs in the original triangle. From there, using Lemma \ref{triangle}, one gets to the simpler case.

$\bullet$ {\it Induction basis}: suppose that the arc distance from $P$ to $Q$ is equal to zero. This means that there is only one arc partially dashed from $P$ to $Q$.

If there is a crossing with both branches dashed, consider the first such crossing encountered when traveling from $P$ according to the orientation, say $X$.
Let $a$ be the part of $D$ that is covered while traveling in the net from the upper branch to the lower branch of $X$. Let $D''$ be the diagram obtained from $D$ by erasing $a$. Let $w'$ be the sum of signs of crossings at which both branches belong to $a$ (including $X$).
Then, as $D$ is descending from the upper to the lower branch of $X$, $H(D)=(xv^{-1})^{w'} H(D'')$ by Lemma \ref{easy_shortening}.
And, by induction on $k$, $H(D'')=(xv^{-1})^{w-w'} H(D')$.

Now, if there are no crossings with both branches dashed, the situation is as in one of two cases shown in Figure \ref{1arc}.

\begin{figure}[ht]
\scalebox{1}{\includegraphics{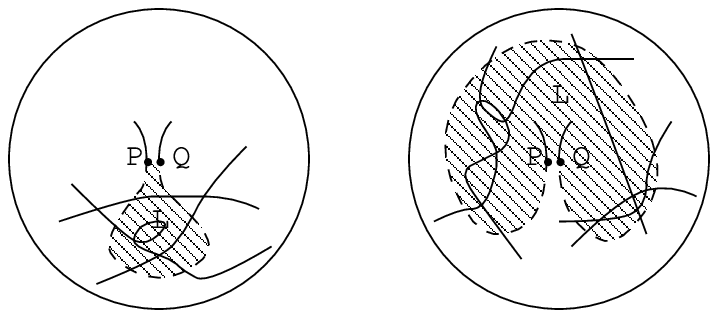}}
\caption{}
\label{1arc}
\end{figure}

In these cases the dashed part is above everything else.
Again, it is possible to order all non-dashed subarcs inside the dashed part $L$. If there are 1-gons inside $L$ they are eliminated as before. Otherwise, subarcs are removed from $L$ using Lemma \ref{removing_2gons}. One gets finally to a situation where there are no non dashed subarcs inside $L$, except for two subarcs containing $P$ and $Q$ in the second case of Figure \ref{1arc}. Then, one checks easily that $H(D)=H(D')$.
\end{proof}
\end{proposition}

\begin{proposition}[Invariance under $\Omega_2'$ move]
$H$ does not change under bad $\Omega_2'$ move involving diagrams with $n$ and $n-2$ crossings.
\begin{proof}
By Remark \ref{nice_remark} and Lemma \ref{independence_reid}, we may suppose that the based diagram involved in the bad move is descending (except possibly at the crossing involved in the move that is not the basepoint).
In the following figures, $w$ stands for the sum of signs of some self-crossings that are removed.

There are two cases to consider. In the first case one branch is oriented downwards and the other is oriented upwards. $H$ is unchanged as it is shown in Figure \ref{bad2}.

\begin{figure}[ht]
\scalebox{1}{\includegraphics{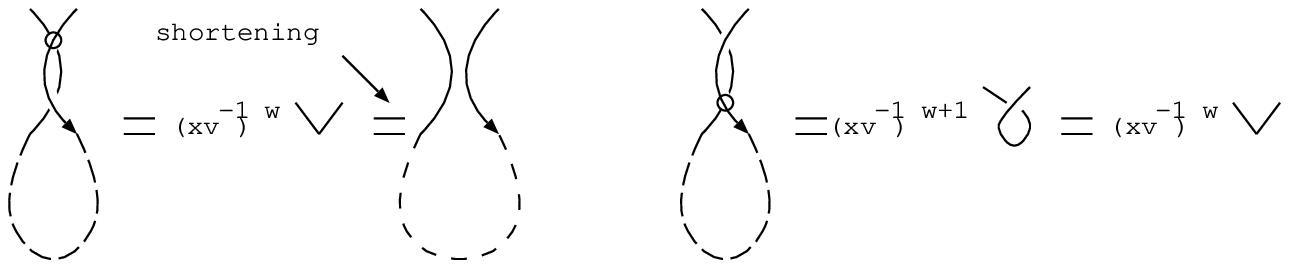}}
\caption{}
\label{bad2}
\end{figure}

In the second case both branches involved in the move are oriented downwards. $H$ is unchanged as it is shown in Figure \ref{bad2bis}.

\begin{figure}[ht]
\scalebox{1}{\includegraphics{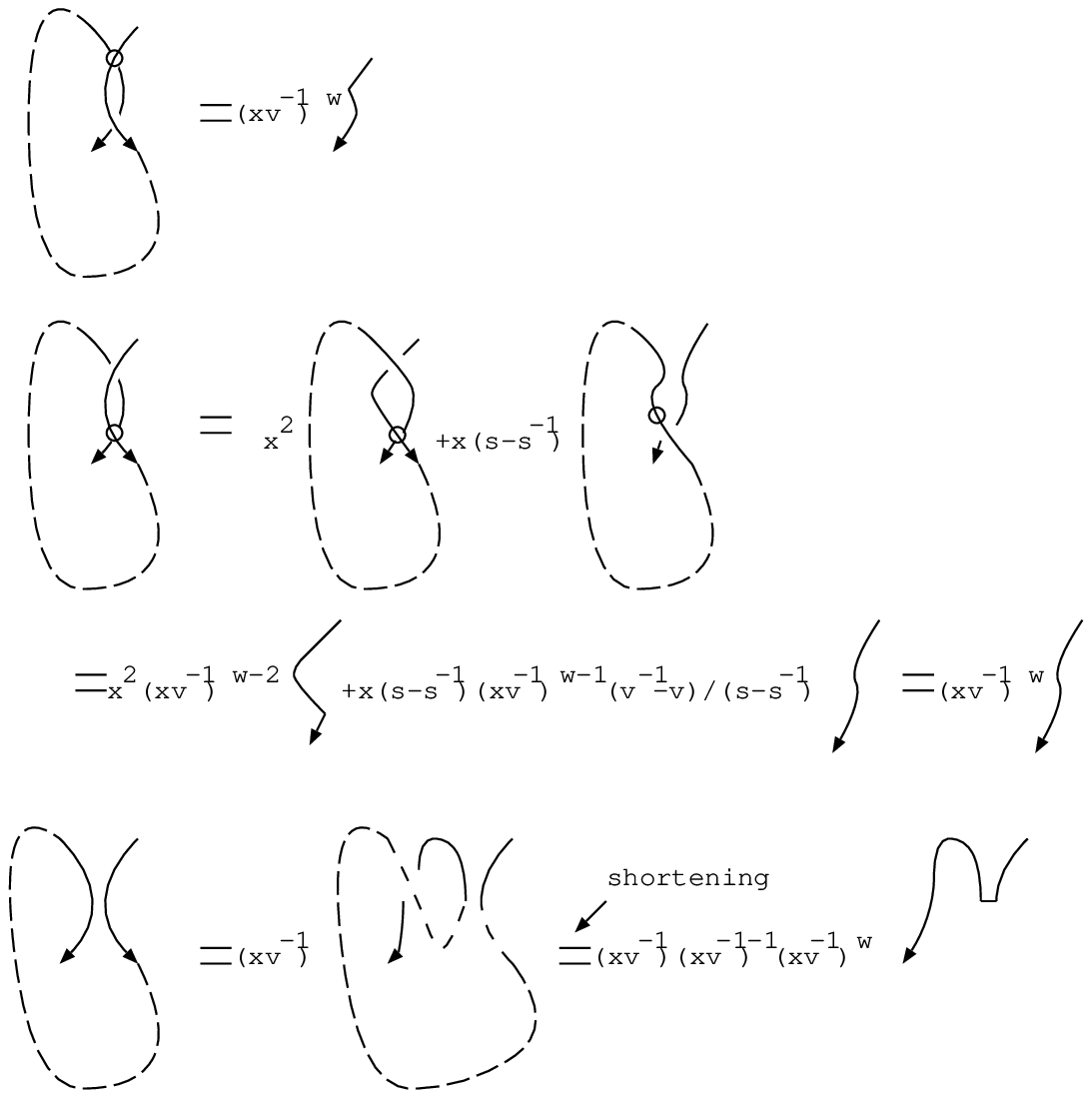}}
\caption{}
\label{bad2bis}
\end{figure}

\end{proof}
\end{proposition}

\begin{proposition}[Invariance under $\Omega_3'$ move]
$H$ does not change under bad $\Omega_3 '$ move involving diagrams with $n$ crossings.
\begin{proof}
By Remark \ref{nice_remark}, one may suppose that the diagrams involved in the move are descending except, possibly, at the crossings appearing in the move.
Furthermore, using Lemma \ref{independence_reid}, one may assume that the diagram before the move is descending.

In the Figures \ref{bad3}, \ref{bad3bis}, \ref{bad3bisbis} and \ref{bad3bisbis_o}, the calculations for $H$ under $\Omega_3 '$ moves are shown. One has also to calculate $H$ under $\Omega_3'$ moves obtained from the moves presented in these figures by performing reflections with respect to vertical lines passing through the basepoints, but the calculations are similar. This gives all possible $\Omega_3'$ moves.

In Figures \ref{bad3bis} and \ref{bad3bisbis} one supposes that the arc distance from the basepoint to the lowest branch is even, whereas in Figure \ref{bad3bisbis_o} it is odd.
$w$, $w_1$ and $w_2$ stand for the sums of signs of some self-crossings that are removed.

\begin{figure}[ht]
\scalebox{1}{\includegraphics{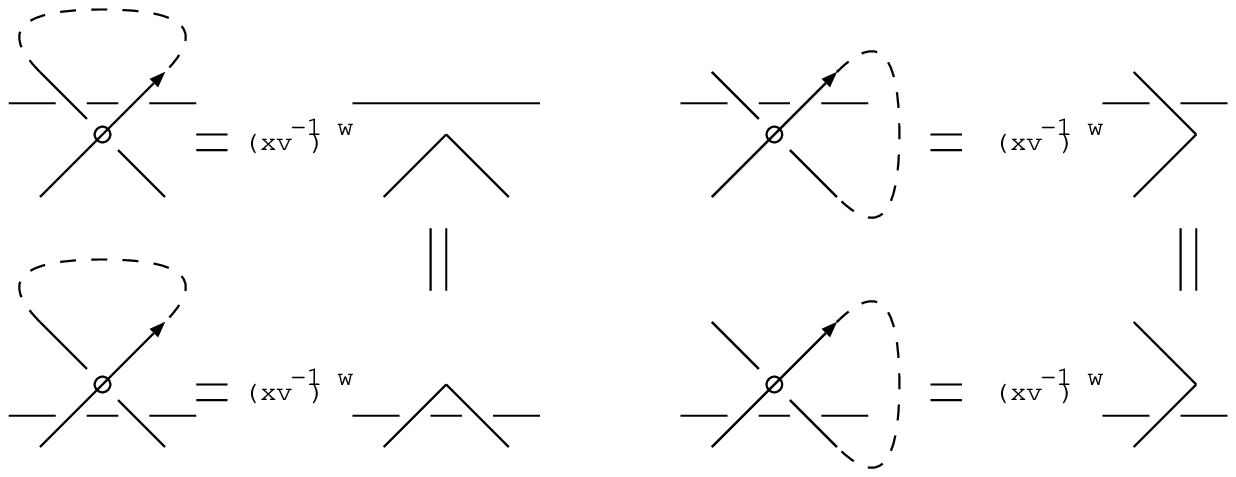}}
\caption{}
\label{bad3}
\end{figure}

\begin{figure}[ht]
\scalebox{1}{\includegraphics{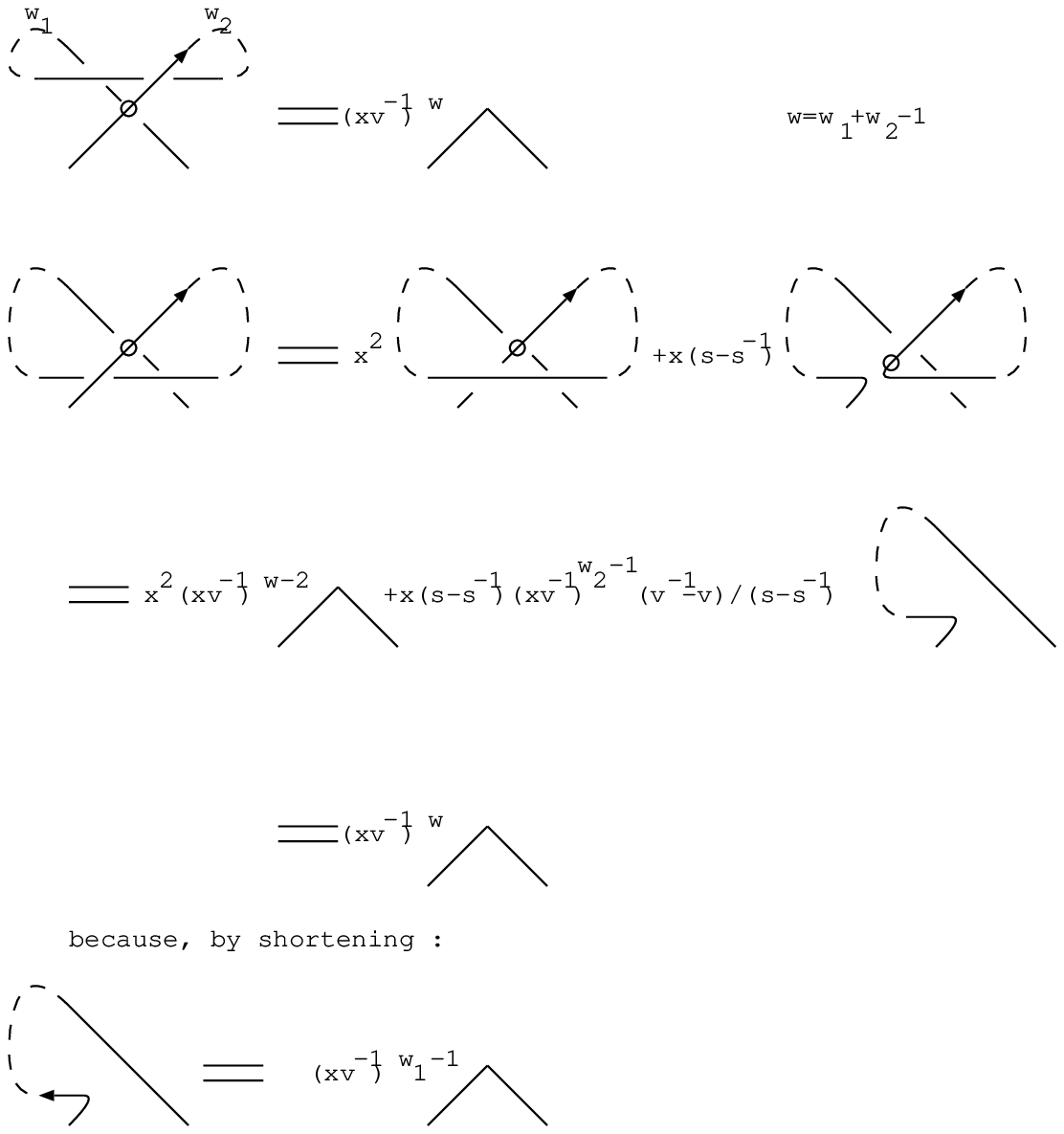}}
\caption{}
\label{bad3bis}
\end{figure}

\begin{figure}[ht]
\scalebox{1}{\includegraphics{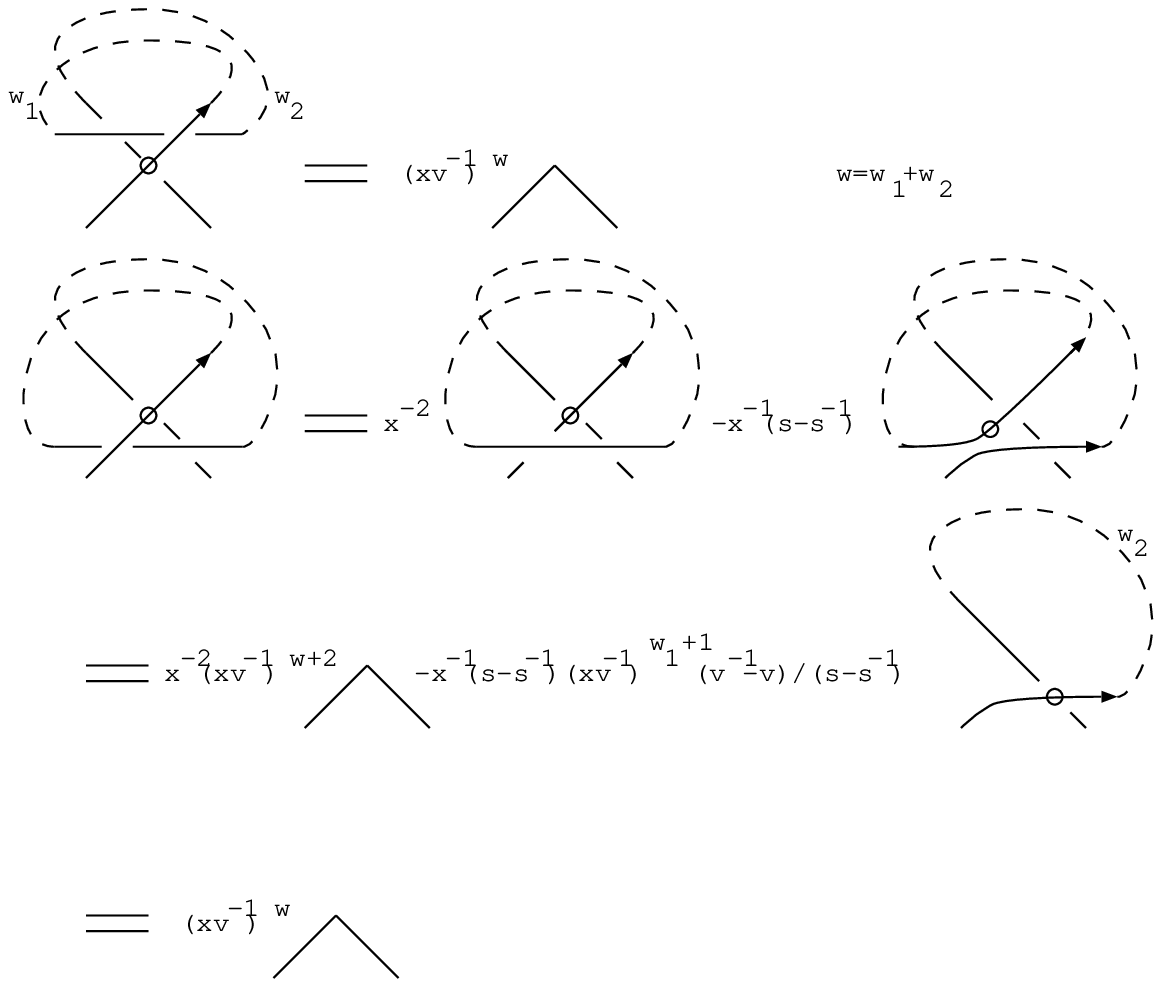}}
\caption{}
\label{bad3bisbis}
\end{figure}

\begin{figure}[ht]
\scalebox{1}{\includegraphics{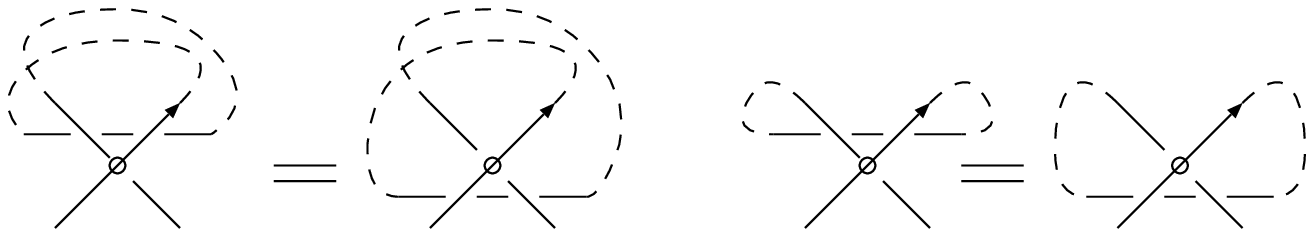}}
\caption{}
\label{bad3bisbis_o}
\end{figure}

\end{proof}
\end{proposition}

\subsection{Bad moves $\Omega_4''$ and $\Omega_5'''$}

A diagram $D$ is said to be {\it ascending} from $P$ to $Q$ if $D'$ is descending from $P$ to $Q$, where $D'$ is the diagram obtained from $D$ by switching every crossing.

For a non simple diagram $\dot{D}$, denote by $w(\dot{D})$ the sum of signs at self-crossings of the 0-homologous component with basepoint; or the sum of signs at self-crossings of the dashed part determined by the basepoint.

\begin{lemma}[shortening of some based diagrams]\label{shortening_n}

Consider based diagrams $\dot{D_1}$ and $\dot{D_2}$ with at most $n$ crossings, presented in Figure \ref{bsp_inf}.

\begin{figure}[ht]
\scalebox{1}{\includegraphics{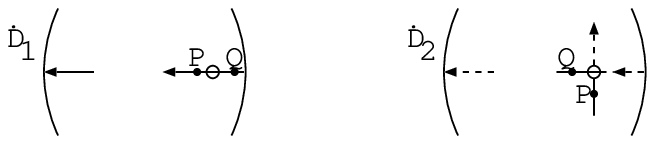}}
\caption{}
\label{bsp_inf}
\end{figure}

Suppose that $\dot{D_1}$ is ascending from $P$ to $Q$.
Let $D_1'$ be the diagram obtained from $\dot{D_1}$ by removing the 0-homologous component with basepoint.
Then $H(\dot{D_1})=\mu (xv^{-1})^{w (\dot{D_1})}H(D_1')$.

Suppose that the dashed part determined by the basepoint of $\dot{D_2}$ is ascending i.e. $\dot{D_2}$ is ascending from $P$ to $Q$.
Let $D_2'$ be the diagram obtained from $\dot{D_2}$ by removing the dashed part determined by the basepoint.
Then $H(\dot{D_2})=(xv^{-1})^{w (\dot{D_2})}H(D_2')$.

\begin{proof}
The proof is similar to the proof of Proposition \ref{shortening}.
It is done by induction on the number of crossings in diagrams $\dot{D_1}$ and $\dot{D_2}$. If this number is 0 in the case of $\dot{D_1}$ or 1 in the case of $\dot{D_2}$ then the lemma follows from the definition of $H$.

Assuming that the lemma is true for diagrams with less then $k$ crossings, it is proven for diagrams with $k$ crossings by induction on $l$: the arc distance from $P$ to $Q$.

It is clear from Figure \ref{bsp_inf} that $l$ is at least equal to 2.
If there are some 1-gons, one uses (HI), induction on $k$ and Lemma \ref{removing_1gons_strong} to eliminate them (see the proof of Proposition \ref{shortening}, {\it Case 1}). {\it Any} 1-gon can be eliminated in this way, because $P$ and $Q$ are not inside any 1-gon.
If there are no 1-gons, one reduces the arc distance from $P$ to $Q$ in the same way as it was done in the proof of Proposition \ref{shortening}, {\it Case 2}, by eliminating the highest arc (the one that is encountered when traveling in the net from $P$ and crossing the line at infinity once), while keeping the part of diagram near $P$ and $Q$ unchanged.

Finally, if $l$ is equal to 2 and there are no 1-gons, one reduces the arc distance from $P$ to $Q$ to 0. Then, using Lemma \ref{moving_arcs}, one reduces the number of crossings involving the component with basepoint to 0 in the case of $\dot{D_1}$; or one reduces the number of crossings involving the dashed part determined by the basepoint to 1 in the case of $\dot{D_2}$. The conclusion of the lemma follows from the definition of $H$.
\end{proof}
\end{lemma}

\begin{proposition}[invariance under $\Omega_4''$ and $\Omega_5'''$ moves]
$H$ does not change under bad moves $\Omega_4''$ and $\Omega_5'''$, involving diagrams with $n$ crossings.
\begin{proof}
Consider first an $\Omega_4''$ move.
Notice that it can be obtained with the crossing of the line at infinity in the net by the basepoint followed by a good $\Omega_4$ move.
It follows from Lemmas \ref{crossing_independence} and \ref{shortening_n} that $H$ does not change when the basepoint crosses the line at infinity in the net (before this crossing the component with basepoint may be assumed descending, and it becomes ascending after the crossing). Also $H$ does not change under good Reidemeister moves. Thus $H$ does not change under an $\Omega_4''$ move.

Notice, that Lemma \ref{shortening_n} is also true if in the diagram $\dot{D_2}$ presented in Figure \ref{bsp_inf}, the vertical branch of the basepoint points downwards instead of upwards.

Consider now an $\Omega_5'''$ move.
It can obtained with a good $\Omega_4$ move and the crossing by the basepoint of the line at infinity, presented in Figure \ref{bad5} (or a similar crossing, for which the vertical branch of the basepoint on the left of Figure \ref{bad5}, points downwards instead of upwards).
It follows from Lemmas \ref{crossing_independence} and \ref{shortening_n} that $H$ does not change when the basepoint crosses the line at infinity in the net (before this crossing the dashed part determined by the basepoint may be assumed descending, and it becomes ascending after the crossing).
Also $H$ does not change under good moves. Thus $H$ does not change under an $\Omega_5'''$ move.

\begin{figure}[ht]
\scalebox{1}{\includegraphics{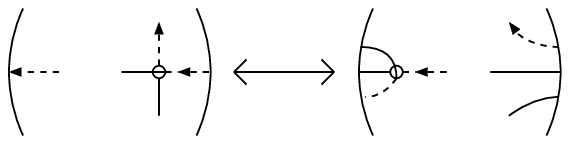}}
\caption{}
\label{bad5}
\end{figure}

\end{proof}
\end{proposition}

\subsection{Basepoints of non simple diagrams}
\begin{proposition}[independence on basepoints for non simple diagrams]\label{independence_nonsimple}
Let $\dot{D}_1$ and $\dot{D}_2$ be two non simple based diagrams with $n$ crossings. Suppose that $\dot{D}_1$ and $\dot{D}_2$ differ only by the position of the basepoint. Then $H(\dot{D}_1)=H(\dot{D}_2)$.
\begin{proof}
First, suppose that there are some 2-gons or 1-gons in $\dot{D}_1$.

Notice that as $H$ does not change under all Reidemeister moves that do not increase the number of crossings beyond $n$, the Lemmas \ref{moving_arcs}, \ref{removing_2gons} and \ref{removing_1gons} can be extended to a situation where the basepoint is allowed to be inside the 2-gon or 1-gon that is removed.
Using Lemma \ref{crossing_independence} and the extended versions of Lemmas \ref{removing_2gons} and \ref{removing_1gons} one may reduce the number of crossings in $\dot{D}_1$ and $\dot{D}_2$. Then $H(\dot{D}_1)=H(\dot{D}_2)$ by IH($n-1$).

It can be seen easily that, if there are no 1-gons and no 2-gons in $\dot{D}_1$, then there has to be an arc with endpoints on the boundary circle of the diagram that are not antipodal. Let $n$ be the number of such arcs in $\dot{D}_1$.

The proof of the proposition is done by induction on $n$. If $n=0$ then there are some 2-gons or 1-gons in $\dot{D}_1$.

If there are $n$ arcs in $\dot{D}_1$ and $\dot{D}_2$, one arc can be removed using Lemma \ref{crossing_independence}, the extended version of Lemma \ref{moving_arcs} and $\Omega_5$ and $\Omega_4$ moves.
Thus, by induction on $n$, $H(\dot{D}_1)=H(\dot{D}_2)$.
\end{proof}
\end{proposition}

This was the final step in the proof of the induction hypothesis IH($n$). Theorem \ref{th_homfly_explicit} follows.

\section{Inductive definition of the Kauffman polynomial $K$}
The proof of Theorem \ref{th_kauffman_explicit} is similar to the proof of Theorem \ref{th_homfly_explicit} above. The differences in the proofs are the same as in the case of links in $\mathbb R^3$. The parts of proofs that are specific to the situation of links in $\mathbb RP^3$ are almost identical for Homfly and Kauffman polynomials.

$K$ is defined in a similar way to $H$ (see section \ref{ind_homfly}). The definition is in fact simpler because in the unoriented case there is a good notion of descending diagram (see \cite{M}).

A {\it standard diagram} of standard unoriented unlink $L_n$ is the diagram presented in Figure \ref{unlink} with the orientations being disregarded.
Let $d=(a+a^{-1}) z^{-1}-1$.

\subsection{Inductive hypothesis IH$(n-1)$}
There is a function $K$ defined on the set of diagrams with at most $(n-1)$ crossings, taking values in $\mathbb Z [a^{\pm 1},z^{\pm 1},y]$ such that:
\begin{enumerate}
\item
$K$ is invariant under those Reidemeister moves that do not increase the number of crossings beyond $n-1$.
\item
$K$ satisfies relations (KI) and (KII).
\item
If $D$ is the standard diagram of standard unoriented unlink $L_m$, $m>0$, with at most $(n-1)$ crossings (i.e. $m(m-1)/2\le n-1$), then $K(D)=y^m$.
Also, $K(\vcenter{\hbox{\epsfig{file=circle_unoriented.eps}}})=d$.
\end{enumerate}

\subsection{Diagrams with no crossings}
As the definition of $K$ uses induction on the number of crossings, $K$ is first defined for diagrams with 0 crossings.

Let $D$ be a diagram with 0 crossings. Let $p$ be the number of its 0-homologous components and $m$ the number of its 1-homologous components ($m$ is 0 or 1). Then, by definition:

$(K_1): K(D)=d ^p y^m$

For convenience $K$ of the empty link is, by definition, equal to 1 (which agrees with $(K_1)$). Note that $K$ satisfies IH(0).

\subsection{Diagrams with $n\ge 1$ crossings}
We assume that the inductive hypothesis IH($n-1$) holds true.

As in the case of the polynomial $H$, the construction of $K$ for diagrams with $n$ crossings is divided into several cases treated in the subsequent subsections.
In each case a diagram $D$ with $n$ crossings is endowed with some extra structure (a {\it directed} basepoint or a couple of basepoints). $D$ together with this structure is denoted by $\dot{D}$. A diagram $\alpha(\dot{D})$ is then defined: it is a diagram obtained from $\dot{D}$ by a series of crossing changes.

The diagram $\alpha(\dot{D})$ has the following property: if $X$ is one of the crossings of $\dot{D}$ that have to be switched to obtain $\alpha(\dot{D})$, and $\dot{D}'$ is the diagram obtained from $\dot{D}$ by switching $X$, then $\alpha(\dot{D})=\alpha(\dot{D}')$.

In the following subsections $K$ is defined on $\alpha(\dot{D})$ for each case (see $(K_3)$ to $(K_5)$).

Suppose that $K$ is already defined on all $\alpha(\dot{D})$. For a based diagram $\dot{D}$, denote by $S(\dot{D})$ the set of crossings of $\dot{D}$ where $\dot{D}$ and $\alpha(\dot{D})$ differ. Let $k$ be the number of elements in $S(\dot{D})$ and $\omega$ a (linear) ordering of $S(\dot{D})$. Denote by $S(\dot{D},\omega)$ the set $S(\dot{D})$ equipped with ordering $\omega$.

We define $K(\dot{D},\omega)$ by induction on $k$. The definition depends on $\omega$.
If $k=0$ then $K$ is already defined.
Otherwise let $\dot{D}'$ be the based diagram obtained from $\dot{D}$ by switching the first crossing in $S(\dot{D},\omega)$; let $D_1''$ and $D_2''$ be the diagrams obtained by smoothing the same crossing in two possible ways. Let $\omega '$ be an ordering of all crossings of $S(\dot{D'})$ induced by $\omega$. Note that there are $k-1$ elements in $S(\dot{D'})$. $K(\dot{D}',\omega ')$ is defined by induction on $k$ and $K(D_1'')$ and $K(D_2'')$ are defined by IH($n-1$).
Now $K(\dot{D},\omega)$ is defined using the relation (KI) on the first crossing in $S(\dot{D},\omega)$ with the help of $K(\dot{D}',\omega ')$, $K(D_1'')$ and $K(D_2'')$. By definition:

$(K_2): K(\dot{D},\omega)=-K(\dot{D}',\omega ')+z (K(D_1'')+K(D_2''))$

\subsection{Simple diagrams}
The definition of simple diagram is the same in oriented and unoriented case (see subsection \ref{simple_diagrams}). The same is true for the definitions of {\it based simple diagram} and {\it primary and secondary basepoints}.

In the case of unoriented links, a based simple diagram $\dot{D}$ is said to be {\it descending} if, making it oriented in any way, it is descending (see bottom of section \ref{simple_diagrams}). It can be easily seen that this definition does not depend on the choice of orientation that is made.
For a based simple diagram $\dot{D}$, let $\alpha(\dot{D})$ be the based diagram obtained from $\dot{D}$ by crossing changes that make it descending.

Let $p$ be the number of 0-homologous components and $m$ the number of 1-homologous components of $\dot{D}$. By definition:

$(K_3): K(\alpha(\dot{D}))=d^p y^m$

\subsection{Non simple diagrams}
In the case of non simple diagrams, a basepoint for an unoriented link is defined as a basepoint for oriented link (see section \ref{nonsimple_diagrams}); moreover, it is endowed with an arrow giving a local orientation at the basepoint. In the case of a basepoint which is a self-crossing of a 1-homologous component the local orientation is given to one of its branches. A basepoint together with an arrow is called {\it directed basepoint}.

The notion of non simple descending diagram in the oriented case depends only on the orientation of the component with basepoint (see section \ref{nonsimple_diagrams}). In the unoriented case, a non simple based diagram $\dot{D}$ is {\it descending} if, endowing the component on which the directed basepoint lies with the orientation given by the arrow of this basepoint, it is descending in the oriented sense. The diagram $\alpha(\dot{D})$ is the based diagram obtained from $\dot{D}$ by the crossing changes which are necessary to make it descending.

Let $\dot{D}$ be a based diagram for which the directed basepoint is on a 0-homologous component.
Let $D'$ be the diagram obtained from $\dot{D}$ by removing this 0-homologous component. Let $w$ be the sum of all signs at all self-crossings of the component with basepoint in $\alpha(\dot{D})$, where this component is oriented arbitrarily. Then, by definition:

$(K_4): K(\alpha({\dot{D}}))=d a^w K(D')$

Let $\dot{D}$ be a based diagram for which the directed basepoint is a self-crossing of a 1-homologous component. Let $D'$ be the diagram obtained from $\dot{D}$ by removing the dashed part determined by the basepoint. Endowing with an arbitrary orientation the component with basepoint, let $w$ be the sum of all signs at all self-crossings of the dashed part determined by the basepoint in $\alpha(\dot{D})$, including the basepoint (which is a self-crossing). Then, by definition:

$(K_5): K(\alpha({\dot{D}}))=a^w K(D')$

\section{Independence of $K$ on choices, invariance under Reidemeister moves}

\subsection{Relations (KI) and (KII)}
\begin{lemma}[independence on ordering]\label{order_independent_k}
Let $\dot{D}$ be a based diagram with $n$ crossings. Let $\omega$ and $\omega'$ be two orderings of the set of crossings of $\dot{D}$ that differ between $\dot{D}$ and $\alpha (\dot{D})$.

Then $K(\dot{D},\omega)=K(\dot{D},\omega')$
\begin{proof}
By induction on the number of crossing differences between $\dot{D}$ and $\alpha (\dot{D})$ it is sufficient to prove that $K$ does not change if one switches the first two crossings according to $\omega$, say $C_1$ and $C_2$.

Denote by $\sigma_i \dot{D}$ the diagram obtained from $\dot{D}$ by switching $C_i$ ($i=1, 2$). Denote by $\mu_i \dot{D}$ and $\nu_i \dot{D}$ the diagrams obtained from $\dot{D}$ by smoothing $C_i$ ($i=1, 2$) in two different ways. Here it does not matter which smoothing is $\mu_i \dot{D}$ and which smoothing is $\nu_i \dot{D}$. First consider the sequence in which $C_1$ is switched before $C_2$:

\begin{eqnarray*}
K(\dot{D},\omega) &=& -K(\sigma_1 \dot{D})+z K(\mu_1 \dot{D})+z K(\nu_1 \dot{D})\\
&=& K(\sigma_2 \sigma_1 \dot{D})-z K(\mu_2 \sigma_1 \dot{D})-z K(\nu_2 \sigma_1 \dot{D})+z K(\mu_1 \dot{D})+z K(\nu_1 \dot{D})
\end{eqnarray*}

And, switching $C_2$ before $C_1$:

\begin{eqnarray*}
K(\dot{D},\omega') &=& -K(\sigma_2 \dot{D})+z K(\mu_2 \dot{D})+z K(\nu_2 \dot{D})\\
&=& K(\sigma_1 \sigma_2 \dot{D})-z K(\mu_1 \sigma_2 \dot{D})-z K(\nu_1 \sigma_2 \dot{D})+z K(\mu_2 \dot{D})+z K(\nu_2 \dot{D})
\end{eqnarray*}

Now:

\begin{eqnarray*}
K(\dot{D},\omega)-K(\dot{D},\omega')=z (-K(\mu_2 \sigma_1 \dot{D})-K(\mu_2 \dot{D})-K(\nu_2 \sigma_1 \dot{D})-K(\nu_2 \dot{D})\\+K(\mu_1 \dot{D})+K(\mu_1 \sigma_2 \dot{D})+K(\nu_1 \dot{D})+K(\nu_1 \sigma_2 \dot{D}))
\end{eqnarray*}

Because of IH($n-1$), one can use (KI) on diagrams with $n-1$ crossings, so:

\begin{eqnarray*}
K(\dot{D},\omega)-K(\dot{D},\omega')=z^2 (-K(\mu_2 \mu_1 \dot{D})-K(\mu_2 \nu_1 \dot{D})-K(\nu_2 \mu_1 \dot{D})-K(\nu_2 \nu_1 \dot{D})\\
+K(\mu_1 \mu_2 \dot{D})+K(\mu_1 \nu_2 \dot{D})+K(\nu_1 \mu_2 \dot{D})+K(\nu_1 \nu_2 \dot{D}))=0
\end{eqnarray*}

\end{proof}
\end{lemma}

Lemma \ref{easy_shortening} can be easily modified to the case of $K$:

\begin{lemma}
Let $D$ be a diagram with at most $n-1$ crossings and let $X$ be a self-crossing of a component $b$ of $D$. Suppose that for some fixed orientation of $b$, the arc distance from the upper branch to the lower branch of $X$ is even and $D$ is descending from the upper branch to the lower branch of $X$.

Let $b'$ be the part of $b$ that is covered if one travels in the net from the upper branch to the lower branch of $X$, according to the fixed orientation of $b$. Let $D'$ be the diagram obtained from $D$ by erasing $b'$. Let $w$ be the sum of signs of self-crossings of $b'$ (including $X$).

Then $K(D)=a^w K(D')$.
\end{lemma}

In the same way as for $H$, from the preceding lemma and Lemma \ref{order_independent_k} follows:

\begin {proposition}[Kauffman relations]
The relation (KII) holds for $K$ in the case when the diagram on the left has $n$ crossings. The relation (KI) holds for $K$ in the case when the two based diagrams on the left have $n$ crossings and these diagrams have the same basepoint(s).
\end{proposition}

\subsection{Basepoints for 0-homologous components}
The Proposition \ref{shortening} for $H$ is a consequence of the inductive hypothesis IH($n-1$) only. It can be easily modified to the case of $K$:

\begin{proposition}[shortening of diagrams]\label{shortening_k}
Let $D$ be a diagram with $k$ crossings. A part of $D$ is shown on the left of Figure \ref{lemma_desc_k}. Suppose that, if the component of $D$ which contains $P$ and $Q$ is oriented in such a way that one can travel in the net from $P$ to $Q$ while covering the dashed part on the left of Figure \ref{lemma_desc_k}, then the arc distance from $P$ to $Q$ is even and $D$ is descending from $P$ to $Q$.

$D'$ is obtained from $D$ by removing the dashed part and joining $P$ and $Q$ with a segment (a part of $D'$ is shown on the right of Figure \ref{lemma_desc_k}). Let $w$ be the sum of signs at all crossings for which both branches are in the dashed part.

If $k\le n-1$ then $K(D)=a^w K(D')$.

\begin{figure}[ht]
\scalebox{1}{\includegraphics{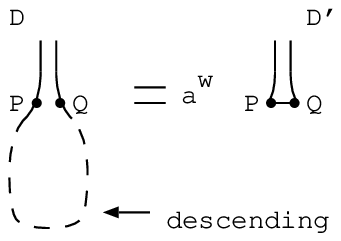}}
\caption{}
\label{lemma_desc_k}
\end{figure}

\end{proposition}

\begin{lemma}[moving the basepoint]
Suppose that $P$ is a directed basepoint lying on a 0-homologous component $b$ of a based diagram $D$ with $n$ crossings and $P$ is on an arc $c$. 
Then $K$ does not change if $P$ is moved on $c$.
\begin{proof}
It is sufficient to prove that $K$ is unchanged if the basepoint passes through a crossing as in Figure \ref{k_basepoint}.

\begin{figure}[ht]
\scalebox{1}{\includegraphics{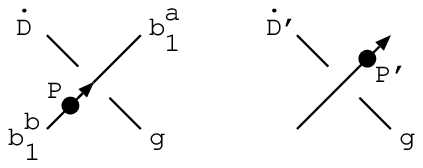}}
\caption{}
\label{k_basepoint}
\end{figure}

Notice that, as relation (KI) holds for diagrams with $n$ crossings, if $K$ does not change for some $D$, when moving the basepoint, then it does not change for any diagram obtained from $D$ by some crossing changes (using IH($n-1$)).
We may therefore suppose that the based diagram $\dot{D}$ on the left of Figure \ref{k_basepoint} is descending.
Starting from $P$ and traveling on the net of $D$ according to the orientation given by the arrow of the basepoint, denote the successive arcs encountered by $b_1, b_2, ..., b_l$. Furthermore denote the part of $b_1$ that comes after $P$ by $b_1^a$ and the remaining part by $b_1^b$. Exactly in the same way as in the proof of Lemma \ref{moving_basepoint} one has the following:

$$b_2\le b_4\le b_6\; ...\le b_1^b\; ...\le b_5\le b_3\le b_1^a$$

where $b_i\le b_j$ means that $b_i$ is below $b_j$.

Now $\dot{D}'$ may be descending or not.
It is not descending if and only if the branch $g$ (see Figure \ref{k_basepoint}) is a part of $b_1^a$, $b_1^b$ or $b_k$ with $k$ odd.
If $\dot{D}'$ is not descending, it becomes descending if one switches the crossing in \ref{k_basepoint}.

Let $\sigma \dot{D}'$, $\eta \dot{D}'$ and $\eta_2 \dot{D}'$ be the diagrams obtained from $\dot{D}'$ by doing respectively the switching at this crossing, the smoothing at this crossing respecting any orientation of $b$, and the other possible smoothing at the crossing. Note that in $\eta \dot{D}'$, $b$ becomes a link with two 0-homologous components $b^1$ and $b^2$ where $b^1$ contains $P'$ and $b^2$ contains $P$ ($P$ and $P'$ can naturally be viewed in $\eta \dot{D}'$). 
Notice that $\eta \dot{D}'$ is descending with respect to $P'$.

Denote by $D_b$ the diagram obtained from $\dot{D}'$ by removing $b$, which is the same as the diagram obtained from $\eta \dot{D}'$ by removing $b^1$ and $b^2$. Let $\epsilon$ be the sign of the crossing in Figure \ref{k_basepoint} and $w$ the sum of signs of crossings for which both branches are in $b$.

Denote by $D_{b^1}$ the diagram obtained from $\eta \dot{D}'$ by removing $b^1$. As $D_{b^1}$ is descending with respect to $P$, one has $K(\eta \dot{D}')=d^2 a^{w-\epsilon} K(D_b)$.

In $\eta_2 \dot{D}'$, there is a single 0-homologous component coming from $b$. It is not descending, as part of it is descending and the other part is ascending.
One uses Proposition \ref{shortening_k} ($\eta_2 \dot{D}'$ has $n-1$ crossings) to eliminate the descending part. $K$ of the remaining ascending part can be computed by changing the direction of the directed basepoint $P$ so that this part becomes descending (by IH($n-1$), $K$ does not depend on the choice of directed basepoint). One gets $K(\eta_2 \dot{D}')=d a^{w-\epsilon} K(D_b)$.

Then:

$K(\dot{D})=d a^w H(D_b)$

$K(\dot{D}')=-K(\sigma \dot{D}')+z(K(\eta \dot{D}')+K(\eta_2 \dot{D}'))=-d a^{w-2\epsilon } K(D_b)+z(d^2 a^{w-\epsilon} K(D_b)+d a^{w-\epsilon} K(D_b))=-d a^{w-2\epsilon} K(D_b)+z d(d a^{w-\epsilon}+a^{w-\epsilon}) K(D_b)=d (-a^{w-2\epsilon}+z (a^{w-\epsilon} (a+a^{-1}) z^{-1}-a^{w-\epsilon}+a^{w-\epsilon})) K(D_b) =d a^w K(D_b)$

as $d=(a+a^{-1}) z^{-1}-1$.

Thus $K(\dot{D})=K(\dot{D}')$.
\end{proof}
\end{lemma}

\subsection{Invariance of $K$ under good and bad Reidemeister moves. Independence on basepoints.}
The invariance of $K$ under good Reidemeister moves is proved similarly to the invariance of $H$. The calculations that have to be done are the same as in the case of Kauffman polynomial for classical links \cite{K}.
From this follows, as it was the case for $H$, the independence of $K$ on basepoints for simple diagrams (this is Proposition \ref{simple_unchanged} modified to $K$).

\begin{figure}[ht]
\scalebox{1}{\includegraphics{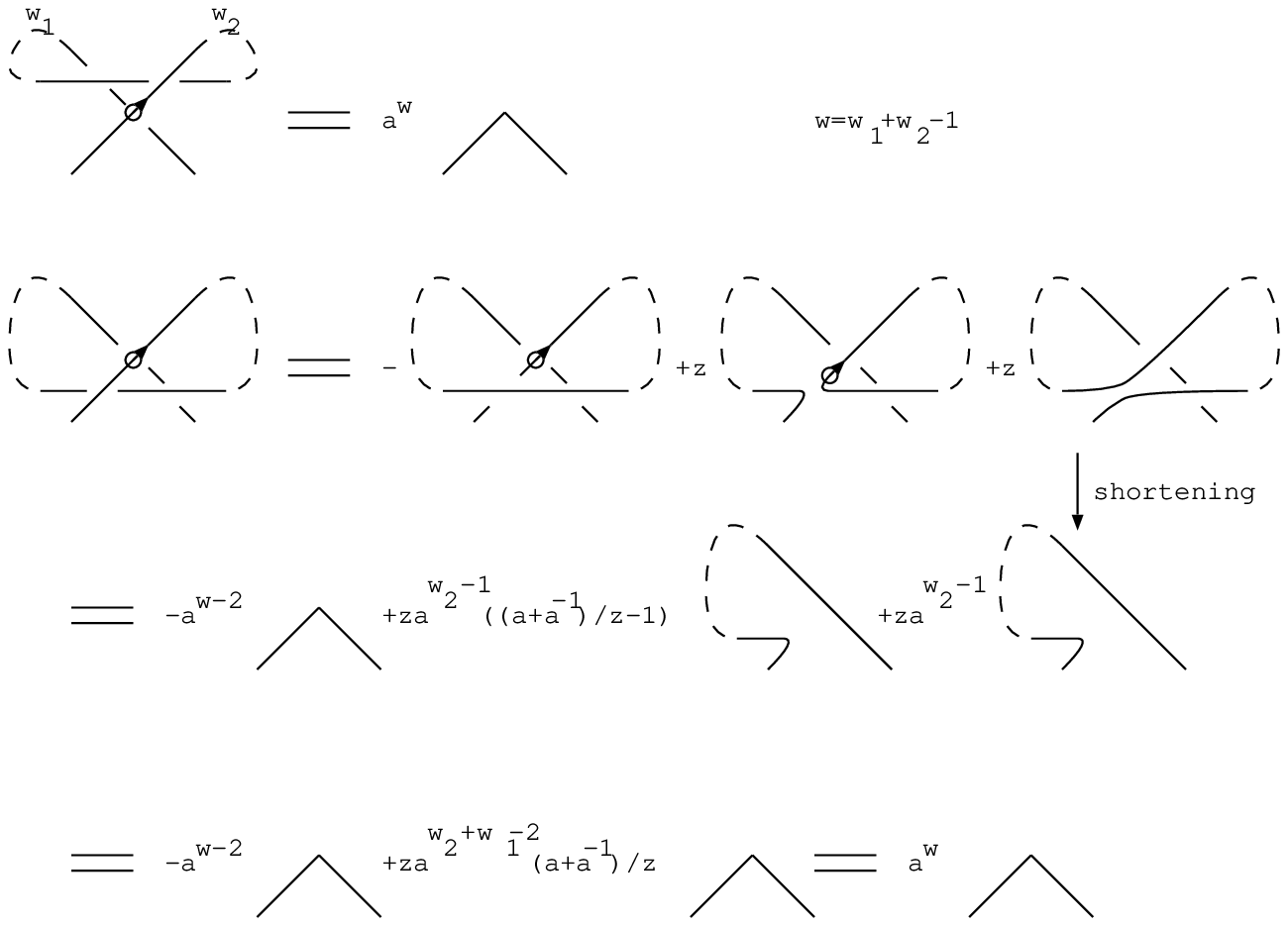}}
\caption{}
\label{bad3_k}
\end{figure}

$K$ is also unchanged under bad Reidemeister moves. An example of calculation is shown in Figure \ref{bad3_k}. In this figure the arc distance from the basepoint to the middle branch involved in the move is even.

The independence of $K$ on directed basepoints for non simple diagrams is proven exactly as it was done for $H$ in Proposition \ref{independence_nonsimple}.

We have established:

\begin{proposition}[invariance under Reidemeister moves, independence on basepoints]
$K$ does not change under any Reidemeister move that involves diagrams with at most $n$ crossings. For diagrams with $n$ crossings, $K$ does not dependent on basepoints.
\end{proposition}

Thus, assuming that $K$ satisfies IH($n-1$), we have shown that it satisfies IH($n$). Theorem \ref{th_kauffman_explicit} follows.

\section{An application: distance from affinity}

The {\it distance from affinity} of a link in $\mathbb RP^3$ is, by definition, the minimum on all its diagrams of the number of times the line at infinity is intersected in the net. For example, a link is affine if and only if its distance from affinity is equal to 0.

The Homfly and Kauffman polynomials can be used to get a lower bound for the distance from affinity of a link:

\begin{proposition}
Let $L$ be a framed oriented link and suppose that $H(L)$ has degree $n$ in $z$. Then the distance from affinity of $L$ is at least equal to $n$.

Let $L$ be a framed unoriented link and suppose that $K(L)$ has degree $n$ in $y$. Then the distance from affinity of $L$ is at least equal to $n$.
\begin{proof}
The proof is the same for $H$ and $K$. Suppose that $L$ is a framed oriented link and that $H(L)$ has degree $n$ in $z$.
Suppose that there is a diagram of $L$ in which the line at infinity is intersected in the net less then $n$ times. Then, computing $H$ on this diagram does not give rise to terms with degree in $z$ greater or equal to $n$, by definition of $H$. But in that case $H(L)$ cannot have degree $n$ in $z$.
\end{proof}
\end{proposition}

\begin{figure}[ht]
\scalebox{1}{\includegraphics{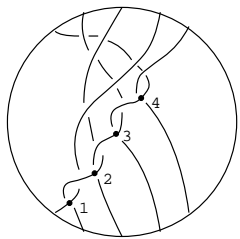}}
\caption{A knot with distance from affinity equal to 5}
\label{knot_5}
\end{figure}

The proposition above can be used to show that, for any $n\in \mathbb N \cup \{ 0\}$, there exists knots with distance from affinity equal to $n$. An example for $n=5$ is shown in Figure \ref{knot_5}. By definition, the distance from affinity of this knot is at most $5$.
To see that it is at least $5$, use the Homfly skein relation (HI) successively on the $4$ crossings marked with a point in this figure, and get $H(L_5)$ with a factor $x^4 (s-s^{-1})^4$ from the smoothings.
For the links coming from crossing changes that appear when using (HI), it can be seen easily that their distance from affinity is at most equal to $3$ so, in $H$, they do not contribute to the term of degree $5$ in $z$.

An interesting question is whether $H$ or $K$ can detect exactly the distance from affinity of any link.

\end{document}